\documentstyle[twoside,psfig,lscape]{article}

\oddsidemargin=\evensidemargin 
\catcode`\@=11 \long\def\@makefntext#1{ \protect\noindent \hbox to
3.2pt {\hskip-.9pt
$^{{\eightrm\@thefnmark}}$\hfil}#1\hfill}       

\def\ps@myheadings{\let\@mkboth\@gobbletwo      
\def\@oddhead{\hbox{}
\rightmark\hfil\eightrm\thepage}
\def\@oddfoot{}\def\@evenhead{\eightrm\thepage\hfil
\leftmark\hbox{}}\def\@evenfoot{}
\def\sectionmark##1{}\def\subsectionmark##1{}}

\catcode`@=11                        
\def\ps@plain{\let\@mkboth\@gobbletwo
     \def\@oddhead{}\def\@oddfoot{\eightrm\hfil\thepage
     \hfil}\def\@evenhead{}\let\@evenfoot\@oddfoot}

\newcounter{sectionc}\newcounter{subsectionc}\newcounter{subsubsectionc}
\renewcommand{\section}[1] {\vspace{12pt}\addtocounter{sectionc}{1}
\setcounter{subsectionc}{0}\setcounter{subsubsectionc}{0}\noindent
    {\tenbf\thesectionc. #1}\par\vspace{5pt}}
\renewcommand{\subsection}[1] {\vspace{12pt}\addtocounter{subsectionc}{1}
    \setcounter{subsubsectionc}{0}\noindent
    {\bf\thesectionc.\thesubsectionc.
    {\kern1pt \bfit #1}}\par\vspace{5pt}}
\renewcommand{\subsubsection}[1] {\vspace{12pt}
    \addtocounter{subsubsectionc}{1}
    \noindent
    {\tenrm\thesectionc.\thesubsectionc.\thesubsubsectionc. {\kern1pt
    \it #1}}\par\vspace{5pt}}

\newcounter{appendixc}
\newcounter{subappendixc}[appendixc]
\newcounter{subsubappendixc}[subappendixc]

\renewcommand{\appendix}[1] {\vspace{12pt}  
    \refstepcounter{appendixc}      
    \setcounter{figure}{0}
    \setcounter{table}{0}
    \setcounter{lemma}{0}
    \setcounter{theorem}{0}
    \setcounter{corollary}{0}
    \setcounter{definition}{0}
    \setcounter{equation}{0}
    \renewcommand{\thefigure}{\Alph{appendixc}.\arabic{figure}}
    \renewcommand{\thetable}{\Alph{appendixc}.\arabic{table}}
    \renewcommand{\theappendixc}{\Alph{appendixc}}
    \renewcommand{\thelemma}{\Alph{appendixc}.\arabic{lemma}}
    \renewcommand{\thetheorem}{\Alph{appendixc}.\arabic{theorem}}
    \renewcommand{\thedefinition}{\Alph{appendixc}.\arabic{definition}}
    \renewcommand{\thecorollary}{\Alph{appendixc}.\arabic{corollary}}
    \renewcommand{\theequation}{\Alph{appendixc}.\arabic{equation}}
    \noindent{\tenbf Appendix \theappendixc #1}\par\vspace{5pt}}

\newcommand{\smalllineskip}{\baselineskip=10pt}

\newcommand{\copyrightheading}[1]
    {\vspace*{-2.5cm}\smalllineskip{\flushleft
    {\footnotesize }\\
    {\footnotesize \copyright\kern2pt }\\
         }}

\def\keywords#1{{
    \centering{\begin{minipage}{4.5in}\footnotesize\baselineskip=10pt
    {\footnotesize\it Keywords}\/: #1
    \end{minipage}}\par}}

\renewenvironment{thebibliography}[1]
    {\frenchspacing
     \ninerm\baselineskip=11pt
     \begin{list}{[\arabic{enumi}]}
    {\usecounter{enumi}\setlength{\parsep}{0pt}
     \setlength{\leftmargin 13.7pt}{\rightmargin 0pt} 
     \setlength{\itemsep}{0pt} \settowidth
    {\labelwidth}{[#1]}\sloppy}}{\end{list}}

\newcounter{itemlistc}
\newcounter{romanlistc}
\newcounter{alphlistc}
\newcounter{arabiclistc}

\newcommand{\fcaption}[1]{
        \refstepcounter{figure}
        \setbox\@tempboxa = \hbox{\footnotesize Fig.~\thefigure. #1}
        \ifdim \wd\@tempboxa > 5in
           {\begin{center}
        \parbox{5in}{\footnotesize\smalllineskip Fig.~\thefigure. #1}
            \end{center}}
        \else
             {\begin{center}
             {\footnotesize Fig.~\thefigure. #1}
              \end{center}}
        \fi}

\newcommand{\tcaption}[1]{
        \refstepcounter{table}
        \setbox\@tempboxa = \hbox{\footnotesize Table~\thetable. #1}
        \ifdim \wd\@tempboxa > 5in
           {\begin{center}
        \parbox{5in}{\footnotesize\smalllineskip Table~\thetable. #1}
            \end{center}}
        \else
             {\begin{center}
             {\footnotesize Table~\thetable. #1}
              \end{center}}
        \fi}

\def\pmb#1{\setbox0=\hbox{#1}
    \kern-.025em\copy0\kern-\wd0
    \kern.05em\copy0\kern-\wd0
    \kern-.025em\raise.0433em\box0}

\def\fnt#1#2{\footnotetext{\kern-.3em
    {$^{\mbox{\scriptsize #1}}$}{#2}}}

\font\tenrm=cmr10  \font\tenbf=cmbx10
\font\bfit=cmbxti10 at 10pt \font\ninerm=cmr9 
 \font\eightrm=cmr8

\newtheorem{theorem}{Theorem}   

\newtheorem{lemma}{Lemma}

\newtheorem{definition}{Definition}

\def\@begintheorem#1#2{\trivlist    
    \item[\hskip\labelsep{\bf #1\ #2.}]}
\def\@opargbegintheorem#1#2#3{\trivlist
    \item[\hskip\labelsep{\bf #1\ #2\ (#3).}]}

        {\setcounter{itemlistc}{0}      
     \begin{list}{$\bullet$}        
    {\usecounter{itemlistc}         
     \leftmargin10pt           
     \setlength{\parsep}{0pt}
     \setlength{\itemsep}{0pt}     
    }}{\end{list}}

    {\setcounter{romanlistc}{0}     
     \begin{list}{$($\roman{romanlistc}$)$} 
    {\usecounter{romanlistc}        
     \leftmargin18pt 
     \setlength{\parsep}{0pt}
     \setlength{\itemsep}{0pt}  
     \settowidth{\labelwidth}{#1}
    }}{\end{list}}

    {\setcounter{enumii}{0}         
     \begin{list}{$($\alph{enumii}$)$}  
    {\usecounter{enumii}            
     \leftmargin18pt        
     \setlength{\parsep}{0pt}
     \setlength{\itemsep}{0pt}  
     \settowidth{\labelwidth}{#1}
    }}{\end{list}}

\def\qed{\hbox{${\vcenter{\vbox{            
   \hrule height 0.4pt\hbox{\vrule width 0.4pt height 6pt
   \kern5pt\vrule width 0.4pt}\hrule height 0.4pt}}}$}}

\def\theequation{\thesectionc.\arabic{equation}}  

\begin{document}

\markboth{Slavik Jablan} {Ascending numbers of knots and their families}

\centerline{\bf UNKNOTTING AND ASCENDING NUMBERS}
\centerline{\bf OF KNOTS AND THEIR FAMILIES}
\bigskip

\centerline{\footnotesize SLAVIK JABLAN}
\medskip

\bigskip

\centerline{\footnotesize\it The Mathematical Institute, Knez
Mihailova 36,}\centerline{\footnotesize\it P.O.Box 367, 11001
Belgrade, Serbia} \centerline{\footnotesize\it sjablan@gmail.com}

\bigskip

\begin{abstract}
Ascending numbers are determined for 64 knots with at most $n=10$ crossings. After proving the theorem about the signature of alternating knot families,
we distinguished all families of knots obtained from generating alternating knots with at most 10 crossings, for which the unknotting number can be confirmed by using the general formulae for signatures. For 11 families of knots general formulae are obtained for their ascending numbers.
\end{abstract}

\keywords{Conway notation, knot family, signature, unknotting number, ascending number.}


\section{Introduction}

About Conway notation of knots the reader can consult the seminal paper by J.~Conway \cite{2}, where this notation is introduced, the paper by A.~Caudron \cite{3}, and books \cite{4,5}. In particular, drawings of all knots up to $n=10$ crossings according to Conway notation, where every knot is represented by a single diagram, are given in the Appendix C of the book "Knots and links" by D.~Rolfsen \cite{4}.

In Sections 2,3 and 4 we compute ascending numbers for 64 knots with at most $n=10$ crossings and determine upper and lower bounds of ascending numbers for all knots up to $n=10$ crossings. For twist knots, i.e., knots of the family $p\,2$ ($p\ge 1$) in the Conway notation, the ascending number is one, and for all other knots $a(K)\ge 2$, i.e., $a(K)\ge max(u(K),2)$. This means, that if there is a diagram $\widetilde{K}$ of a $K$ with $a(\widetilde{K})=u(K)$, then $a(K)=a(\widetilde{K})=u(K)$. Except for several knots, the unknotting numbers of knots with at most $n=10$ crossings are known, and they are given in "Tables of knot invariants" by C.~Livingston and J.C.~Cha \cite{6}. Bridge numbers of knots with $n\le 10$ crossings are given in the same tables, but they are not useful for our purpose because for all knots with $n\le 10$ the bridge number is 2 or 3. In order to improve upper bound given by the inequality (1) following from minimal crossing numbers, we computed ascending numbers of all minimal diagrams. As an additional improvement, for some knots we obtained upper bounds from ascending numbers of some of their non-minimal diagrams. For all computations we used the program "LinKnot" \cite{5}.

In Section 5 we prove the theorem on signature enabling computation of general formulae for the signature of alternating knot families given by their Conway symbols. These general formulae enabled us to recognize the families of knots obtained from alternating generating knots with at most $n=10$ crossings for which unknotting numbers are determined by signatures computed in Section 6.

In Section 7 we consider some families of knots with ascending numbers that coincides with the unknotting number.

The ascending number of a link $L$ is described in the paper "Ascending number of knots and links" by M.~Ozawa \cite{1}. In our paper we restrict the consideration
of ascending numbers to knots, so we repeat the definitions of basic terms from \cite{1}. A knot diagram is {\it based} if a base point (different from the crossing
points) is specified on the diagram, and {\it oriented} if an orientation is assigned to it. Let $K$ be a knot and $\widetilde{K}$ be a based oriented diagram of $K$.
The {\it descending diagram} of $\widetilde{K}$, denoted by $d(\widetilde{K})$, is obtained as follows: beginning at the basepoint of $\widetilde{K}$ and proceeding in
the direction specified by the orientation, change the crossings as necessary so that each crossing is first encountered as an over-crossing. Note that $d(\widetilde{K})$ is the diagram of a trivial knot.

\begin{definition}
Let $K$ be a knot and let $\widetilde{K}$ be a based oriented diagram of $K$. The {\it ascending number} of $K$ is defined as the number of different crossings between
$\widetilde{K}$ and $d(\widetilde{K})$ and denoted by $a(\widetilde{K})$. The ascending number of $K$ is defined as the minimum number of $a(\widetilde{K})$ over all based oriented knot diagrams $\widetilde{K}$ of $K$, and denoted by $a(K)$ \cite{1}.
\end{definition}

Among theorems proved in \cite{1}, we relate four of them giving upper and lower bounds for ascending numbers of knots:
\begin{enumerate}
\item for a non-trivial knot $K$, we have
$$a(K)\le \lfloor{{c(K)-1}\over 2}\rfloor \eqno (1)$$
\noindent where $c(K)$ is the minimum crossing number of $K$, and $\lfloor x\rfloor $ integer part of $x$;
\item for every non-trivial knot $K$, we have
$$a(K)\ge u(K)$$
\noindent where $u(K)$ is the unknotting number of $K$;
\item the ascending number of a knot $K$ is one {\it iff} $K$ is a twist knot;
\item for a knot $K$, we have
$$a(K)\ge b(K)-1$$
\noindent where $b(K)$ is the bridge number of $K$.
\end{enumerate}

\section{Ascending numbers of knots up to 8 crossings}

Ascending numbers of knots up to $n=8$ crossings are given in the tables from paper \cite{1} and illustrated by the corresponding based oriented knot diagrams
giving the minimal ascending number, where the knot diagrams, which are the same as minimal crossing diagrams, are omitted. Among knots with the minimal diagram giving the ascending number we recognized two more knots: $7_6=2\,2\,1\,2$ and $8_{12}=2\,2\,2\,2$, illustrated in Fig. 1. For knots $8_{16}=.2.2\,0$ and $8_{17}=.2.2$ we succeeded to find their non-minimal diagrams with diagram ascending number equal 2, so $a(8_{16})=2$ and $a(8_{17})=2$. In the corresponding tables every knot is given in classical Conway notation {\it Con} \cite{4}, followed by unknotting number $u$ \cite{6}, upper bound for ascending number $a_d$ (obtained mostly from minimal diagrams), and ascending number $a$. For knots with unknown ascending numbers a sequence is given, beginning with lower bounds and ending with the best known upper bound (e.g., $[2,3]$). For knots up to $n=8$ crossings, the computation of ascending numbers corresponding to all minimal diagrams gives no improvement of the upper bound obtained from the crossing number, but for many knots with $n=9$ or $n=10$ crossings it results in the upper bound equals 3, instead of the upper bound 4 obtained from the crossing number.

\begin{figure}[th]
\centerline{\psfig{file=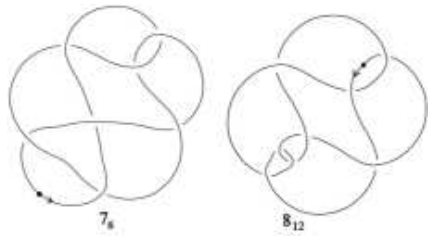,width=2.20in}}
\vspace*{8pt}
\caption{(a) Knot $7_6$; (b) knot $8_{12}$.\label{fig1}}
\end{figure}

\begin{figure}[th]
\centerline{\psfig{file=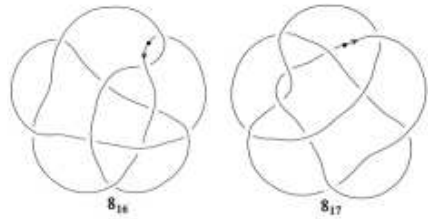,width=2.20in}}
\vspace*{8pt}
\caption{(a) Knot $8_{16}$; (b) knot $8_{17}$.\label{fig2}}
\end{figure}

\footnotesize

\noindent \begin{tabular}{|c|c|c|c|c||c|c|c|c|c|} \hline

$K$ & $Con$ & $u$ & $a_d$ & $a$ & $K$ & $Con$ & $u$ & $a_d$ & $a$   \\ \hline
$3_1$ & $3$ & 1 & 1 & 1 & $7_1$ & $7$ & 3 & 3 & 3  \\ \hline
$4_1$ & $2\,2$ & 1 & 1 & 1 & $7_2$ & $5\,2$ & 1 & 3 & 1  \\ \hline
$5_1$ & $5$ & 2 & 2 & 2 & $7_3$ & $4\,3$ & 2 & 3 &  2 \\ \hline
$5_2$ & $3\,2$ & 1 & 2 & 1 & $7_4$ & $3\,1\,3$ & 2 & 3 & 2  \\ \hline
$6_1$ & $4\,2$ & 1 & 2 & 1 & $7_5$ & $3\,2\,2$ & 2 & 3 & 2  \\ \hline
$6_2$ & $3\,1\,2$ & 1 & 2 & 2 & $7_6$ & $2\,2\,1\,2$ & 1 & 2 & 2  \\ \hline
$6_3$ & $2\,1\,1\,2$ & 1 & 2 & 2 & $7_7$ & $2\,1\,1\,1\,2$ & 1 & 2 & 2  \\ \hline

\end{tabular}

\normalsize

\bigskip

\scriptsize

\noindent \begin{tabular}{|c|c|c|c|c||c|c|c|c|c|} \hline

$K$ & $Con$ & $u$ & $a_d$ & $a$ & $K$ & $Con$ & $u$ & $a_d$ & $a$   \\ \hline
$8_1$ & $6\,2$ & 1 & 3 & 1 & $8_{12}$ & $2\,2\,2\,2$ & 2 & 2 & 2  \\ \hline
$8_2$ & $5\,1\,2$ & 2 & 3 & $[2,3]$ & $8_{13}$ & $3\,1\,1\,1\,2$ & 1 & 3 & 2  \\ \hline
$8_3$ & $4\,4$ & 2 & 3 & 2 & $8_{14}$ & $2\,2\,1\,1\,2$ & 1 & 3 & 2  \\ \hline
$8_4$ & $4\,1\,3$ & 2 & 3 & 2 & $8_{15}$ & $2\,1,2\,1,2$ & 2 & 3 & 2  \\ \hline
$8_5$ & $3,3,2$ & 2 & 3 & $[2,3]$ & $8_{16}$ & $.2.2\,0$ & 2 & 3 & 2  \\ \hline
$8_6$ & $3\,3\,2$ & 2 & 3 & 2 & $8_{17}$ & $.2.2$ & 1 & 3 & 2  \\ \hline
$8_7$ & $4\,1\,1\,2$ & 1 & 3 & $[2,3]$ & $8_{18}$ & $8^*$ & 2 & 2 & 2  \\ \hline
$8_8$ & $2\,3\,1\,2$ & 2 & 3 & 2 & $8_{19}$ & $3,3,-2$ & 3 & 3 & 3  \\ \hline
$8_9$ & $3\,1\,1\,3$ & 1 & 3 & $[2,3]$ & $8_{20}$ & $3,2\,1,-2$ & 1 & 2 & 2  \\ \hline
$8_{10}$ & $2\,1,3,2$ & 2 & 3 & $[2,3]$ & $8_{21}$ & $2\,1,2\,1-2$ & 1 & 2 & 2  \\ \hline
$8_{11}$ & $3\,2\,1\,2$ & 1 & 3 & 2 &  &  &  &  &   \\ \hline

\end{tabular}

\normalsize

\bigskip

\section{Ascending numbers of knots with 9 crossings}

According to \cite{1}, for knots with $n=9$ crossings ascending numbers are known only for six knots: $a(9_3)=3$, $a(9_4)=2$, $a(9_6)=3$,
$a(9_7)=2$, $a(9_{47})=2$, and $a(9_{48}=2$, and they were determined by M.~Okuda. For knots $9_{47}$ and $9_{48}$ they can be determined from their minimal diagrams.
Hence, for 23 new non-trivial knots\footnote{For knots $9_1$ and $9_2$ is trivial to conclude that $a(9_1)=4$, and $a(9_2)=1$.} with $n=9$ crossings we obtained their ascending numbers. Based oriented diagrams corresponding to these knots are illustrated in Figs. 3-10. All these alternating knots with $n=9$ crossings are given by their non-minimal based oriented diagrams giving their ascending numbers. For all remaining knots with $n=9$ crossings, except for the knot $9_{40}$, by computing diagram ascending numbers for all minimal crossing diagrams (or for some non-alternating diagrams in the case of knots $9_{29}$ and $9_{39}$) we succeeded to reduce the set of possible values of the ascending number to $[2,3]$ (meaning, 2 or 3).

\bigskip

\scriptsize

\noindent \begin{tabular}{|c|c|c|c|c||c|c|c|c|c|} \hline

$K$ & $Con$ & $u$ & $a_d$ & $a$ & $K$ & $Con$ & $u$ & $a_d$ & $a$   \\ \hline
$9_1$ & $9$ & 4 & 4 & 4 & $9_{26}$ & $3\,1\,1\,1\,1\,2$ & 1 & 3 & $[2,3]$  \\ \hline
$9_2$ & $7\,2$ & 1 & 4 & 1 & $9_{27}$ & $2\,1\,2\,1\,1\,2$ & 1 & 3 & $[2,3]$  \\ \hline
$9_3$ & $6\,3$ & 3 & 4 & 3 & $9_{28}$ & $2\,1,2\,1,2+$ & 1 & 3 & $[2,3]$  \\ \hline
$9_4$ & $5\,4$ & 2 & 4 & 2 & $9_{29}$ & $.2.2\,0.2$ & 2 & 4 & $[2,3]$  \\ \hline
$9_5$ & $5\,1\,3$ & 2 & 4 & 2 & $9_{30}$ & $2\,1\,1,2\,1,2$ & 1 & 3 & $[2,3]$  \\ \hline
$9_6$ & $5\,2\,2$ & 3 & 4 & 3 & $9_{31}$ & $2\,1\,1\,1\,1\,1\,2$ & 2 & 3 & $[2,3]$  \\ \hline
$9_7$ & $3\,4\,2$ & 2 & 4 & 2 & $9_{32}$ & $.2\,1.2\,0$ & 2 & 3 & $[2,3]$  \\ \hline
$9_8$ & $2\,4\,1\,2$ & 2 & 3 & 2 & $9_{33}$ & $.2\,1.2$ & 1 & 3 & $[2,3]$  \\ \hline
$9_9$ & $4\,2\,3$ & 3 & 4 & 3 & $9_{34}$ & $8^*2\,0$ & 1 & 3 & 2  \\ \hline
$9_{10}$ & $3\,3\,3$ & 3 & 4 & 3 & $9_{35}$ & $3,3,3$ & 3 & 4 & 3  \\ \hline
$9_{11}$ & $4\,1\,2\,2$ & 2 & 3 & $[2,3]$ & $9_{36}$ & $2\,2,3,2$ & 2 & 3 & $[2,3]$  \\ \hline
$9_{12}$ & $4\,2\,1\,2$ & 1 & 3 & 2 & $9_{37}$ & $3,2\,1,2\,1$ & 2 & 3 & 2  \\ \hline
$9_{13}$ & $3\,2\,1\,3$ & 3 & 4 & 3 & $9_{38}$ & $.2.2.2$ & 3 & 4 & 3  \\ \hline
$9_{14}$ & $4\,1\,1\,1\,2$ & 1 & 3 & 2 & $9_{39}$ & $2:2:2\,0$ & 1 & 4 & $[2,3]$  \\ \hline
$9_{15}$ & $2\,3\,2\,2$ & 2 & 3 & 2 & $9_{40}$ & $9^*$ & 2 & 4 & $[2,3,4]$  \\ \hline
$9_{16}$ & $3,3,2+$ & 3 & 4 & 3 & $9_{41}$ & $2\,0:2\,0:2\,0$ & 2  & 3 & $[2,3]$  \\ \hline
$9_{17}$ & $2\,1\,3\,1\,2$ & 2 & 3 & $[2,3]$ & $9_{42}$ & $2\,2,3,-2$ & 1 & 2 & 2  \\ \hline
$9_{18}$ & $3\,2\,2\,2$ & 2 & 4 & 2 & $9_{43}$ & $2\,1\,1,3,-2$ & 2 & 3 & $[2,3]$  \\ \hline
$9_{19}$ & $2\,3\,1\,1\,2$ & 1 & 3 & 2 & $9_{44}$ & $2\,2,2\,1,-2$ & 1 & 2 & 2  \\ \hline
$9_{20}$ & $3\,1\,2\,1\,2$ & 2 & 3 & $[2,3]$ & $9_{45}$ & $2\,1\,1,2\,1,-2$ & 1 & 2 & 2  \\ \hline
$9_{21}$ & $3\,1\,1\,2\,2$ & 1 & 3 & 2 & $9_{46}$ & $3,3,-3$ & 2 & 2 & 2  \\ \hline
$9_{22}$ & $2\,1\,1,3,2$ & 1 & 3 & $[2,3]$ & $9_{47}$ & $8^*-2\,0$ & 2 & 2 & 2  \\ \hline
$9_{23}$ & $2\,2\,1\,2\,2$ & 2 & 4 & 2 & $9_{48}$ & $2\,1,2\,1,-3$ & 2 & 2 & 2  \\ \hline
$9_{24}$ & $3,2\,1,2$ & 1 & 3 & $[2,3]$ & $9_{49}$ & $-2\,0:-2\,0:-2\,0$ & 3 & 3 & 3  \\ \hline
$9_{25}$ & $2\,2,2\,1,2$ & 2 & 3 & 2 &  & &  &  &   \\ \hline
\end{tabular}

\bigskip

\normalsize

\begin{figure}[th]
\centerline{\psfig{file=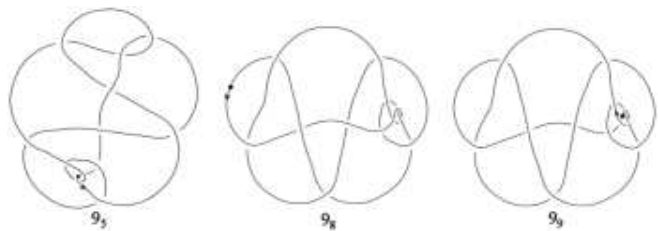,width=3.40in}}
\vspace*{8pt}
\caption{(a) Knot $9_5$; (b) knot $9_8$; (c) knot $9_9$.\label{fig3}}
\end{figure}

\begin{figure}[th]
\centerline{\psfig{file=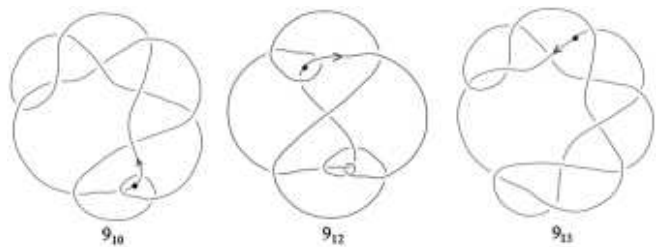,width=3.40in}}
\vspace*{8pt}
\caption{(a) Knot $9_{10}$; (b) knot $9_{12}$; (c) knot $9_{13}$.\label{fig4}}
\end{figure}

\begin{figure}[th]
\centerline{\psfig{file=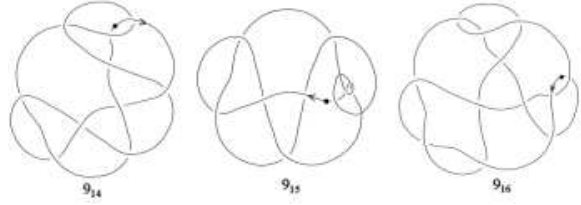,width=3.0in}}
\vspace*{8pt}
\caption{(a) Knot $9_{14}$; (b) knot $9_{15}$; (c) knot $9_{16}$.\label{fig5}}
\end{figure}

\begin{figure}[th]
\centerline{\psfig{file=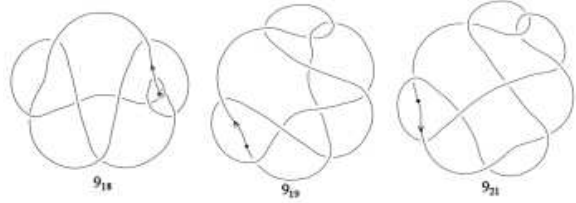,width=3.0in}}
\vspace*{8pt}
\caption{(a) Knot $9_{18}$; (b) knot $9_{19}$; (c) knot $9_{21}$.\label{fig6}}
\end{figure}

\begin{figure}[th]
\centerline{\psfig{file=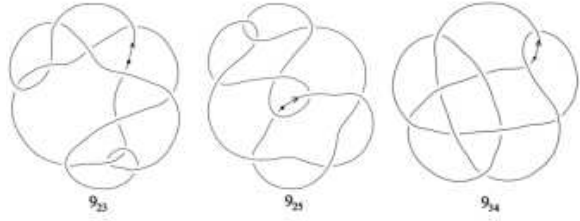,width=3.0in}}
\vspace*{8pt}
\caption{(a) Knot $9_{23}$; (b) knot $9_{25}$; (c) knot $9_{34}$.\label{fig7}}
\end{figure}

\begin{figure}[th]
\centerline{\psfig{file=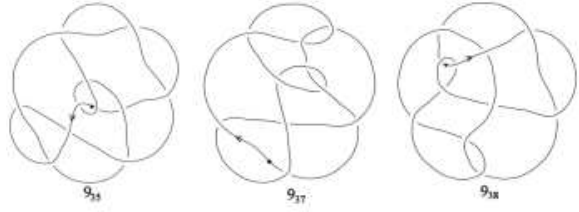,width=3.0in}}
\vspace*{8pt}
\caption{(a) Knot $9_{35}$; (b) knot $9_{37}$; (c) knot $9_{38}$.\label{fig8}}
\end{figure}

\begin{figure}[th]
\centerline{\psfig{file=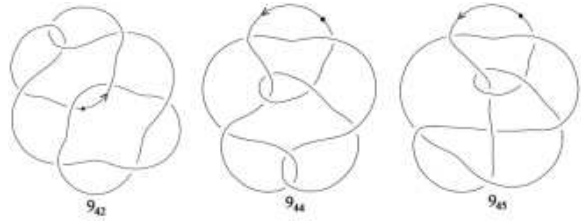,width=3.0in}}
\vspace*{8pt}
\caption{(a) Knot $9_{42}$; (b) knot $9_{44}$; (c) knot $9_{45}$.\label{fig9}}
\end{figure}

\begin{figure}[th]
\centerline{\psfig{file=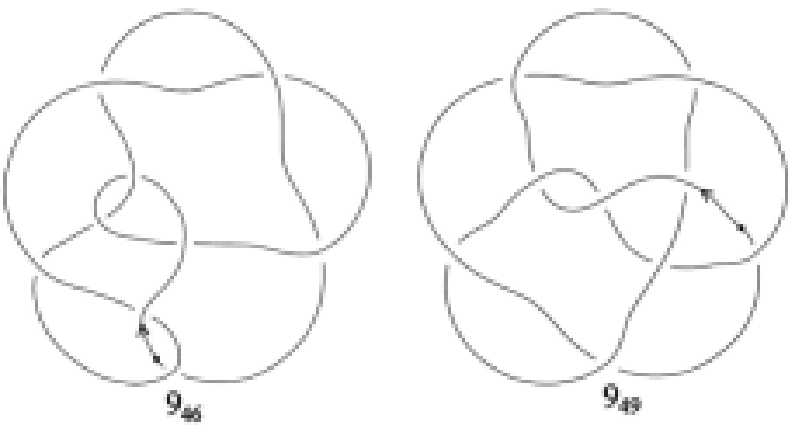,width=2.0in}}
\vspace*{8pt}
\caption{(a) Knot $9_{46}$; (b) knot $9_{49}$.\label{fig10}}
\end{figure}

\section{Ascending numbers of knots with 10 crossings}

For knots with $n=10$ crossings, none of ascending numbers (except for the twist knot $10_1=8\,2$ with $a(10_1)=1$) were known. In this paper we computed ascending numbers for 39 knots with $n=10$ crossings. For some of remaining knots, by using all minimal or some non-minimal diagrams we succeeded to improve upper and lower bound for
ascending numbers to the set $[2,3]$.

Because among 10-crossing knots there are some with unknown unknotting number ($[2,3]$, meaning 2 or 3), the corresponding bounds for ascending number are denoted by $(2,3)$ instead of $[2,3]$, and $(2,3,4)$ instead of $[2,3,4]$. If in any of these cases unknotting number is equal to its lower bound, this will be a counterexample to the Bernhard-Jablan Conjecture \cite{5,7,8}.

Based oriented diagrams of 38 knots with $n=10$ crossings for which we succeeded to compute their ascending numbers are illustrated in Figs. 11-23.

\scriptsize

\noindent \begin{tabular}{|c|c|c|c|c||c|c|c|c|c|} \hline

$K$ & $Con$ & $u$ & $a_d$ & $a$ & $K$ & $Con$ & $u$ & $a_d$ & $a$   \\ \hline
$10_1$ & $8\,2$ & 1 & 4 & 1 & $10_{84}$ & $.2\,2.2$ & 1 & 4 & $[2,3,4]$  \\ \hline
$10_2$ & $7\,1\,2$ & 3 & 4 & $[3,4]$ & $10_{85}$ & $.4.2\,0$ & 2 & 4 & $[2,3]$  \\ \hline
$10_3$ & $6\,4$ & 2 & 4 & 2 & $10_{86}$ & $.3\,1.2\,0$ & 2 & 4 & $[2,3]$  \\ \hline
$10_4$ & $6\,1\,3$ & 2 & 4 & 2 & $10_{87}$ & $.2\,2.2\,0$ & 2 & 4 & $[2,3,4]$  \\ \hline
$10_5$ & $6\,1\,1\,2$ & 2 & 4 & $[2,3,4]$ & $10_{88}$ & $.2\,1.2\,1$ & 1 & 3 & $[2,3]$  \\ \hline
$10_6$ & $5\,3\,2$ & 3 & 4 & 3 & $10_{89}$& $2\,1.2\,1\,0$ & 2 & 3 & $[2,3]$  \\ \hline
$10_7$ & $5\,2\,1\,2$ & 1 & 4 & 2 & $10_{90}$ & $.3.2.2$ & 2 & 4 & $[2,3,4]$  \\ \hline
$10_8$& $5\,1\,4$ & 2 & 3 & $[2,3]$ & $10_{91}$ & $.3.2.2\,0$ & 1 & 4 & $[2,3]$  \\ \hline
$10_9$ & $5\,1\,1\,3$ & 2 & 3 & $[2,3]$ & $10_{92}$ & $.2\,1.2.2\,0$ & 2 & 4 & $[2,3,4]$  \\ \hline
$10_{10}$ & $5\,1\,1\,1\,2$ & 1 & 4 & 2 & $10_{93}$ & $.3.2\,0.2$ & 2 & 4 & $[2,3]$  \\ \hline
$10_{11}$ & $4\,3\,3$ & $[2,3]$ & 4 & $(2,3)$ & $10_{94}$ & $.3\,0.2.2$ & 2 & 4 & $[2,3]$  \\ \hline
$10_{12}$ & $4\,3\,1\,2$ & 2 & 4 & $[2,3]$ & $10_{95}$ & $.2\,1\,0.2.2$ & 1 & 4 & $[2,3,4]$  \\ \hline
$10_{13}$ & $4\,2\,2\,2$ & 2 & 3 & 2 & $10_{96}$ & $.2.2\,1.2$ & 2 & 4 & $[2,3,4]$  \\ \hline
$10_{14}$ & $4\,2\,1\,1\,2$ & 2 & 4 & $[2,3]$ & $10_{97}$ & $.2.2\,1\,0.2$ & 2 & 4 & $[2,3,4]$  \\ \hline
$10_{15}$ & $4\,1\,3\,2$ & 2 & 4 & $[2,3]$ & $10_{98}$ & $.2.2.2.2\,0$ & 2 & 4 & $[2,3,4]$  \\ \hline
$10_{16}$ & $4\,1\,2\,3$ & 2 & 4 & $[2,3]$ & $10_{99}$& $.2.2.2\,0.2\,0$ & 2 & 4 & $[2,3,4]$  \\ \hline
$10_{17}$ & $4\,1\,1\,4$ & 1 & 4 & $[2,3,4]$ & $10_{100}$ & $3:2:2$ & $[2,3]$ & 4 & $(2,3)$  \\ \hline
$10_{18}$& $4\,1\,1\,2\,2$ & 1 & 4 & 2 & $10_{101}$ & $2\,1:2:2$ & 3 & 4 & $[3,4]$  \\ \hline
$10_{19}$ & $4\,1\,1\,1\,3$ & 2 & 4 & $[2,3]$ & $10_{102}$ & $3:2:2\,0$ & 1 & 4 & $[2,3,4]$  \\ \hline
$10_{20}$ & $3\,5\,2$ & 2 & 4 & 2 & $10_{103}$ & $3\,0:2:2$ & 3 & 4 & $3$  \\ \hline
$10_{21}$ & $3\,4\,1\,2$ & 2 & 4 & $[2,3]$ & $10_{104}$ & $3:2\,0:2\,0$ & 1 & 4 & $[2,3,4]$  \\ \hline
$10_{22}$ & $3\,3\,1\,3$ & 2 & 4 & $[2,3]$ & $10_{105}$ & $2\,1:2\,0:2\,0$ & 2 & 4 & $[2,3,4]$  \\ \hline
$10_{23}$ & $3\,3\,1\,1\,2$ & 1 & 4 & $[2,3]$ & $10_{106}$ & $3\,0:2:2\,0$ & 2 & 4 & $[2,3]$  \\ \hline
$10_{24}$ & $3\,2\,3\,2$ & 2 & 4 & 2 & $10_{107}$ & $2\,1\,0:2:2\,0$ & 1 & 4 & $[2,3,4]$  \\ \hline
$10_{25}$ & $3\,2\,2\,1\,2$ & 2 & 4 & $[2,3]$ & $10_{108}$ & $3\,0:2\,0:2\,0$ & 2 & 4 & $[2,3]$  \\ \hline
$10_{26}$ & $3\,2\,1\,1\,3$ & 1 & 4 & $[2,3]$ & $10_{109}$& $2.2.2.2$ & 2 & 4 & $[2,3,4]$  \\ \hline
$10_{27}$ & $3\,2\,1\,1\,1\,2$ & 1 & 4 & $[2,3]$ & $10_{110}$ & $2.2.2.2\,0$ & 2 & 4 & $[2,3,4]$  \\ \hline
$10_{28}$& $3\,1\,3\,1\,2$ & 2 & 4 & $[2,3]$ & $10_{111}$ & $2.2.2\,0.2$ & 2 & 4 & $[2,3,4]$  \\ \hline
$10_{29}$ & $3\,1\,2\,2\,2$ & 2 & 3 & $[2,3]$ & $10_{112}$ & $8^*3$ & 2 & 3 & $[2,3]$  \\ \hline
$10_{30}$ & $3\,1\,2\,1\,1\,2$ & 1 & 4 & $[2,3]$ & $10_{113}$ & $8^*2\,1$ & 1 & 3 & $[2,3]$  \\ \hline
$10_{31}$ & $3\,1\,1\,3\,2$ & 1 & 4 & 2 & $10_{114}$ & $8^*3\,0$ & 1 & 3 & $[2,3]$  \\ \hline
$10_{32}$ & $3\,1\,1\,1\,2\,2$ & 2 & 3 & $[2,3]$ & $10_{115}$ & $8^*2\,0.2\,0$ & 2 & 4 & $[2,3,4]$  \\ \hline
$10_{33}$ & $3\,1\,1\,1\,1\,3$ & 1 & 4 & $[2,3]$ & $10_{116}$ & $8^*2:2$ & 2 & 4 & $[2,3]$  \\ \hline
$10_{34}$ & $2\,5\,1\,2$ & 2 & 4 & 2 & $10_{117}$ & $8^*2:2\,0$ & 2 & 4 & $[2,3]$  \\ \hline
$10_{35}$ & $2\,4\,2\,2$ & 2 & 4 & 2 & $10_{118}$ & $8^*2:.2$ & 1 & 4 & $[2,3,4]$  \\ \hline
$10_{36}$ & $2\,4\,1\,1\,2$ & 2 & 4 & 2 & $10_{119}$& $8^*2:.2\,0$ & 1 & 4 & $[2,3,4]$  \\ \hline
$10_{37}$ & $2\,3\,3\,2$ & 2 & 4 & 2 & $10_{120}$ & $8^*2\,0::2\,0$ & 3 & 4 & $[3,4]$  \\ \hline
$10_{38}$& $2\,3\,1\,2\,2$ & 2 & 4 & 2 & $10_{121}$ & $9^*2\,0$ & 2 & 4 & $[2,3,4]$  \\ \hline
$10_{39}$ & $2\,2\,3\,1\,2$ & 2 & 4 & $[2,3]$ & $10_{122}$ & $9^*.2\,0$ & 2 & 3 & $[2,3]$  \\ \hline
$10_{40}$ & $2\,2\,2\,1\,1\,2$ & 2 & 4 & $[2,3]$ & $10_{123}$ & $10^*$ & 2 & 3 & $[2,3]$  \\ \hline

\end{tabular}

\bigskip

\noindent \begin{tabular}{|c|c|c|c|c||c|c|c|c|c|} \hline

$K$ & $Con$ & $u$ & $a_d$ & $a$ & $K$ & $Con$ & $u$ & $a_d$ & $a$   \\ \hline
$10_{41}$ & $2\,2\,1\,2\,1\,2$ & 2 & 3 & $[2,3]$ & $10_{124}$ & $5,3,-2$ & 4 & 4 & 4  \\ \hline
$10_{42}$ & $2\,2\,1\,1\,1\,1\,2$ & 1 & 3 & $[2,3]$ & $10_{125}$ & $5,2\,1,-2$ & 2 & 3 & $[2,3]$  \\ \hline
$10_{43}$ & $2\,1\,2\,2\,1\,2$ & 2 & 3 & $[2,3]$ & $10_{126}$ & $4\,1,3,-2$ & 2 & 3 & $[2,3]$  \\ \hline
$10_{44}$ & $2\,1\,2\,1\,1\,1\,2$ & 1 & 3 & $[2,3]$ & $10_{127}$ & $4\,1,2\,1,-2$ & 2 & 3 & $[2,3]$  \\ \hline
$10_{45}$ & $2\,1\,1\,1\,1\,1\,1\,2$ & 2 & 3 & $[2,3]$ & $10_{128}$ & $3\,2,3,-2$ & 3 & 3 & 3  \\ \hline
$10_{46}$ & $5,3,2$ & 3 & 4 & $[3,4]$ & $10_{129}$& $3\,2,2\,1,-2$ & 1 & 3 & $[2,3]$  \\ \hline
$10_{47}$ & $5,2\,1,2$ & $[2,3]$ & 4 & $(2,3,4)$ & $10_{130}$ & $3\,1\,1,3,-2$ & 1 & 3 & $[2,3]$  \\ \hline
$10_{48}$& $4\,1,3,2$ & 3 & 4 & $[2,3,4]$ & $10_{131}$ & $3\,1\,1,2\,1,-2$ & 1 & 3 & $[2,3]$  \\ \hline
$10_{49}$ & $4\,1,2\,1,2$ & 3 & 4 & 3 & $10_{132}$ & $2\,3,3,-2$ & 1 & 3 & 2  \\ \hline
$10_{50}$ & $3\,2,3,2$ & 2 & 4 & $[2,3,4]$ & $10_{133}$ & $2\,3,2\,1-2$ & 1 & 3 & 2  \\ \hline
$10_{51}$ & $3\,2,2\,1,2$ & $[2,3]$ & 4 & $(2,3,4)$ & $10_{134}$ & $2\,2\,1,3,-2$ & 3 & 3 & 3  \\ \hline
$10_{52}$ & $3\,1\,1,3,2$ & 2 & 4 & $[2,3]$ & $10_{135}$ & $2\,2\,1,2\,1,-2$ & 2 & 3 & 2  \\ \hline
$10_{53}$ & $3\,1\,1,2\,1,2$ & 3 & 4 & 3 & $10_{136}$ & $2\,2,2\,2-2$ & 1 & 2 & 2  \\ \hline
$10_{54}$ & $2\,3,3,2$ & $[2,3]$ & 4 & $(2,3,4)$ & $10_{137}$ & $2\,2,2\,1\,1,-2$ & 1 & 2 & 2  \\ \hline
$10_{55}$ & $2\,3,2\,1,2$ & 2 & 4 & 2 & $10_{138}$ & $2\,1\,1,2\,1\,1,-2$ & 2 & 3 & $[2,3]$  \\ \hline
$10_{56}$ & $2\,2\,1,3,2$ & 2 & 4 & $[2,3]$ & $10_{139}$& $4,3,-2\,1$ & 4 & 4 & 4  \\ \hline
$10_{57}$ & $2\,2\,1,2\,1,2$ & 2 & 4 & $[2,3]$ & $10_{140}$ & $4,3,-3$ & 2 & 3 & $[2,3]$  \\ \hline
$10_{58}$& $2\,2,2\,2,2$ & 2 & 3 & $[2,3]$ & $10_{141}$ & $4,2\,1,-3$ & 1 & 3 & $[2,3]$  \\ \hline
$10_{59}$ & $2\,2,2\,1\,1,2$ & 1 & 3 & $[2,3]$ & $10_{142}$ & $3\,1,3,-2\,1$ & 3 & 4 & 3  \\ \hline
$10_{60}$ & $2\,1\,1,2\,1\,1,2$ & 1 & 3 & $[2,3]$ & $10_{143}$ & $3\,1,3,-3$ & 1 & 3 & $[2,3]$  \\ \hline
$10_{61}$ & $4,3,3$ & $[2,3]$ & 4 & $(2,3,4)$ & $10_{144}$ & $3\,1,2\,1,-3$ & 2 & 3 & $[2,3]$  \\ \hline
$10_{62}$ & $4,3,2\,1$ & 2 & 4 & $[2,3,4]$ & $10_{145}$ & $2\,2,3,-2\,1$ & 2 & 3 & 2  \\ \hline
$10_{63}$ & $4,2\,12\,1$ & 2 & 4 & $[2,3]$ & $10_{146}$ & $2\,2,2\,1,-3$ & 1 & 2 & 2  \\ \hline
$10_{64}$ & $3\,1,3,3$ & 2 & 4 & $[2,3,4]$ & $10_{147}$ & $2\,1\,1,3,-3$ & 1 & 2 & 2  \\ \hline
$10_{65}$ & $3\,1,3,2\,1$ & 2 & 4 & $[2,3,4]$ & $10_{148}$ & $(3,2)\,(3,-2)$ & 2 & 3 & $[2,3]$  \\ \hline
$10_{66}$ & $3\,1,2\,1,2\,1$ & 3 & 4 & $[3,4]$ & $10_{149}$& $(3,2) \,(2\,1,-2)$ & 2 & 3 & $[2,3]$  \\ \hline
$10_{67}$ & $2\,2,3,2\,1$ & 2 & 4 & $[2,3,4]$ & $10_{150}$ & $(2\,1,2)\,(3,-2)$ & 2 & 3 & $[2,3]$  \\ \hline
$10_{68}$& $2\,1\,1,3,3$ & 2 & 4 & $[2,3,4]$ & $10_{151}$ & $(2\,1,2)\,(2,1,-2)$ & 2 & 3 & $[2,3]$  \\ \hline
$10_{69}$ & $2\,1\,1,2\,1,2\,1$ & 2 & 3 & $[2,3]$ & $10_{152}$ & $(3,2)\,-(3,2)$ & 4 & 4 &  4 \\ \hline
$10_{70}$ & $2\,2,3,2+$ & 2 & 3 & $[2,3]$ & $10_{153}$ & $(3,2)\,-(2\,1,2)$ & 2 & 4 & $[2,3,4]$  \\ \hline
$10_{71}$ & $2\,2,2\,1,2+$ & 1 & 3 & $[2,3]$ & $10_{154}$ & $(2\,1,2)\,-(2\,1,2)$ & 3 & 4 & $[3,4]$  \\ \hline
$10_{72}$ & $2\,1\,1,3,2+$ & 2 & 4 & $[2,3,4]$ & $10_{155}$ & $-3:2:2$ & 2 & 3 & $[2,3]$  \\ \hline
$10_{73}$ & $2\,1\,1,2\,1,2+$ & 1 & 3 & $[2,3]$ & $10_{156}$ & $-3:2:2\,0$ & 1 & 3 & $[2,3]$  \\ \hline
$10_{74}$ & $3,3,2\,1+$ & 2 & 4 & $[2,3,4]$ & $10_{157}$ & $-3:2\,0:2\,0$ & 2 & 3 & $[2,3]$  \\ \hline
$10_{75}$ & $2\,1,2\,1,2\,1+$ & 2 & 3 & $[2,3]$ & $10_{158}$ & $-3\,0:2:2$ & 2 & 3 & $[2,3]$  \\ \hline
$10_{76}$ & $3,3,2++$ & $[2,3]$ & 4 & $(2,3)$ & $10_{159}$& $-3\,0:2:2\,0$ & 1 & 2 & 2  \\ \hline
$10_{77}$ & $3,2\,1,2++$ & $[2,3]$ & 4 & $(2,3)$ & $10_{160}$ & $-3\,0:2\,0:2\,0$ & 2 & 2 & 2  \\ \hline
$10_{78}$& $2\,1,2\,1,2++$ & 2 & 3 & $[2,3]$ & $10_{161}$ & $3:-2\,0:-2\,0$ & 3 & 3 & 3  \\ \hline
$10_{79}$ & $(3,2)\,(3,2)$ & $[2,3]$ & 4 & $(2,3,4)$ & $10_{162}$ & $-3\,0:-2\,0:-2\,0$ & 2 & 3 & $[2,3]$  \\ \hline
$10_{80}$ & $(3,2)\,(2\,1,2)$ & 3 & 4 & 3 & $10_{163}$ & $8^*-3\,0$ & 2 & 2 & 2  \\ \hline
$10_{81}$ & $(2\,1,2)\,(2\,1,2)$ & 2 & 4 & $[2,3]$ & $10_{164}$ & $8^*2:-2\,0$ & 1 & 3 & $[2,3]$  \\ \hline
$10_{82}$ & $.4.2$ & 1 & 4 & $[2,3]$ & $10_{165}$ & $8^*2:.-2\,0$ & 2 & 3 & $[2,3]$  \\ \hline
$10_{83}$ & $.3\,1.2$ & 2 & 4 & $[2,3,4]$ & $$ & $$ &  &  &   \\ \hline

\end{tabular}

\bigskip

\normalsize

\begin{figure}[th]
\centerline{\psfig{file=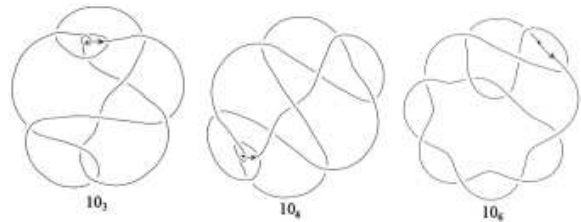,width=3.0in}}
\vspace*{8pt}
\caption{(a) Knot $10_3$; (b) knot $10_4$; (c) knot $10_6$.\label{fig11}}
\end{figure}

\begin{figure}[th]
\centerline{\psfig{file=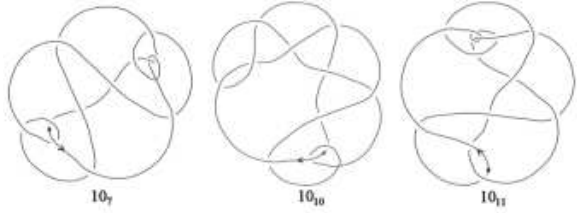,width=3.0in}}
\vspace*{8pt}
\caption{(a) Knot $10_7$; (b) knot $10_{10}$; (c) knot $10_{11}$.\label{fig12}}
\end{figure}

\begin{figure}[th]
\centerline{\psfig{file=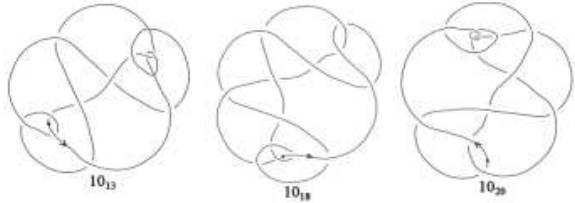,width=3.0in}}
\vspace*{8pt}
\caption{(a) Knot $10_{13}$; (b) knot $10_{18}$; (c) knot $10_{20}$.\label{fig13}}
\end{figure}

\begin{figure}[th]
\centerline{\psfig{file=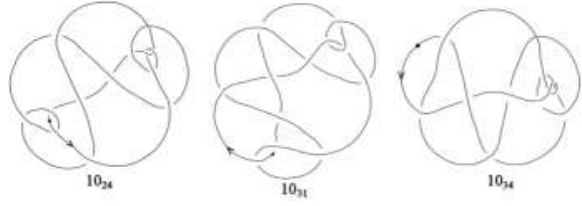,width=3.0in}}
\vspace*{8pt}
\caption{(a) Knot $10_{24}$; (b) knot $10_{31}$; (c) knot $10_{34}$.\label{fig14}}
\end{figure}

\begin{figure}[th]
\centerline{\psfig{file=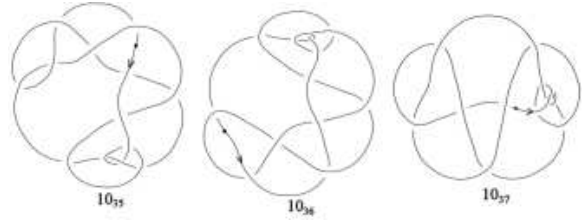,width=3.0in}}
\vspace*{8pt}
\caption{(a) Knot $10_{35}$; (b) knot $10_{36}$; (c) knot $10_{37}$.\label{fig15}}
\end{figure}

\begin{figure}[th]
\centerline{\psfig{file=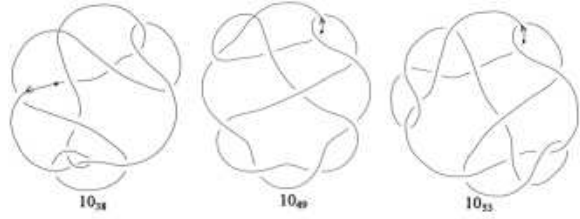,width=3.0in}}
\vspace*{8pt}
\caption{(a) Knot $10_{38}$; (b) knot $10_{49}$; (c) knot $10_{53}$.\label{fig16}}
\end{figure}

\begin{figure}[th]
\centerline{\psfig{file=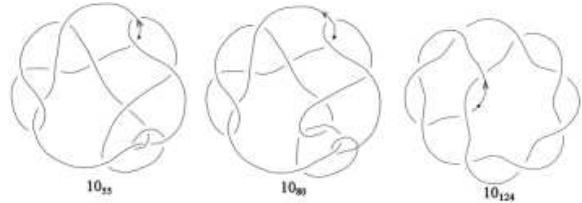,width=3.0in}}
\vspace*{8pt}
\caption{(a) Knot $10_{55}$; (b) knot $10_{80}$; (c) knot $10_{124}$.\label{fig17}}
\end{figure}

\begin{figure}[th]
\centerline{\psfig{file=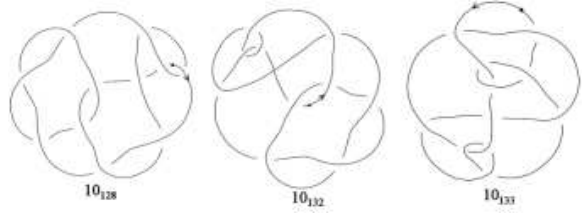,width=3.0in}}
\vspace*{8pt}
\caption{(a) Knot $10_{128}$; (b) knot $10_{132}$; (c) knot $10_{133}$.\label{fig18}}
\end{figure}

\begin{figure}[th]
\centerline{\psfig{file=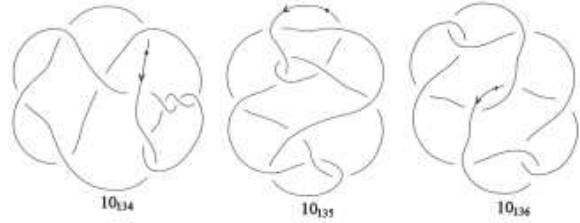,width=3.0in}}
\vspace*{8pt}
\caption{(a) Knot $10_{134}$; (b) knot $10_{135}$; (c) knot $10_{136}$.\label{fig19}}
\end{figure}

\begin{figure}[th]
\centerline{\psfig{file=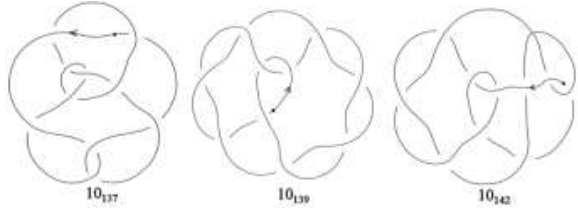,width=3.0in}}
\vspace*{8pt}
\caption{(a) Knot $10_{137}$; (b) knot $10_{139}$; (c) knot $10_{142}$.\label{fig20}}
\end{figure}

\begin{figure}[th]
\centerline{\psfig{file=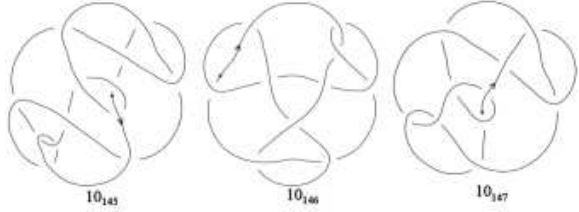,width=3.0in}}
\vspace*{8pt}
\caption{(a) Knot $10_{145}$; (b) knot $10_{146}$; (c) knot $10_{147}$.\label{fig21}}
\end{figure}

\begin{figure}[th]
\centerline{\psfig{file=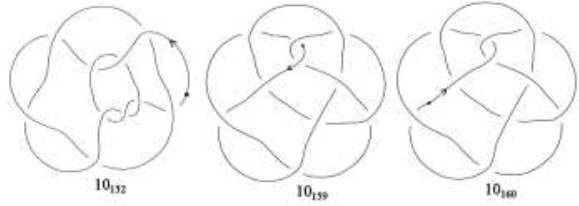,width=3.0in}}
\vspace*{8pt}
\caption{(a) Knot $10_{152}$; (b) knot $10_{159}$; (c) knot $10_{160}$.\label{fig22}}
\end{figure}

\begin{figure}[th]
\centerline{\psfig{file=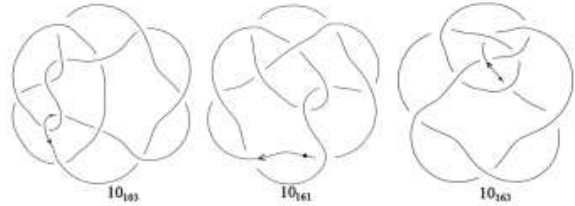,width=3.0in}}
\vspace*{8pt}
\caption{(a) Knot $10_{103}$; (b) knot $10_{161}$; (b) knot $10_{163}$. \label{fig23}}
\end{figure}

\section{Signature and alternating knot families}

\begin{definition}
Let $S$ denote the set of numbers in the unreduced\footnote{The
Conway notation is called {\it unreduced} if $1$'s denoting elementary
tangles in vertices are not omitted in symbols of polyhedral links.}
Conway symbol $C(L)$ of a link $L$. Given $C(L)$ and an arbitrary
(non-empty) subset $\tilde S=\{a_1, a_2, \ldots, a_m \}$ of $S$, the
family $F_{\tilde S}(L)$ of knots or links derived from $L$ is
constructed by substituting each $a_i \in \tilde S$, $a_i \neq 1$ in
$C(L)$ by $sgn(a) (|a|+n)$, for $n \in {N^+}$.
\end{definition}

For even integers $n \geq 0$ this construction preserves the number of
components, i.e., we obtain (sub)families of links with the same number
of components. If all parameters in a Conway symbol of a knot or link  are 1,2, or 3, such a link is
called {\it generating}.

K. Murasugi \cite{9} defined {\it signature} $\sigma _K$ of a knot $K$ as the signature of the matrix $S_K+{S_K}^T$, where ${S_K}^T$ is the transposed matrix of $S_K$, and $S_K$ is the Seifert matrix of the knot $K$.

For alternating knots, signature can be computed by using a combinatorial formula derived by P.~Traczyk \cite{10}. We will use this formula, proved by J.~Przytycki,  in the following form, taken from  \cite{11}, Theorem 7.8, Part (2):

\begin{theorem}
If $D$ is a reduced alternating diagram of an oriented knot, then $$\sigma_D=-{1\over 2}w+{1\over 2}(W-B)=-{1\over 2}w+{1\over 2}(|D_{s+}|-|D_{s-}|),$$ where $w$ is the writhe of $D$, $W$ is the number of white regions in the checkerboard coloring of $D$, which is for alternating minimal diagrams equal to the number of cycles $|D_{s+}|$ in the state $s+$, and $B$ is the number of black regions in the checkerboard coloring of $D$ equal to the number of the cycles $|D_{s-}|$ in the state $s-$.
\end{theorem}

Introducing orientation of a knot, every $n$-twist (chain of digons) becomes
{\it parallel} or {\it anti-parallel}. For signs of crossings and checkerboard coloring we use the conventions shown in Fig. 24.

\begin{figure}[th]
\centerline{\psfig{file=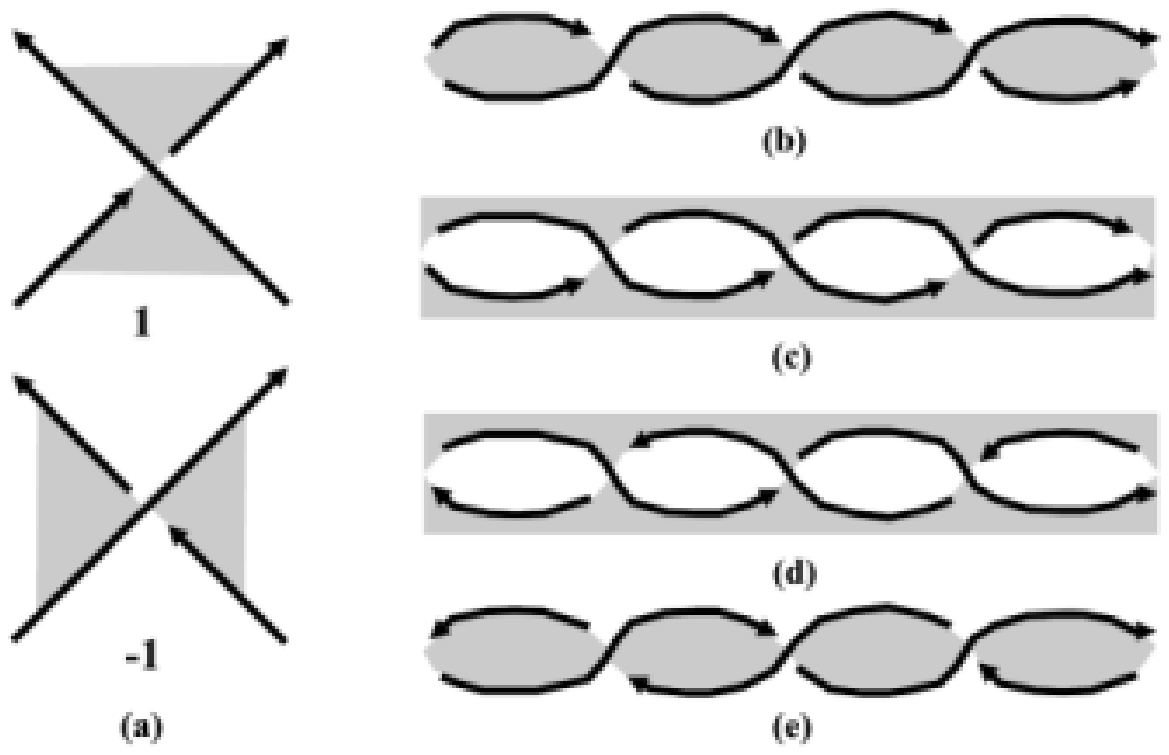,width=1.80in}} \vspace*{8pt}
\caption{(a) Positive crossing and negative crossing (b)
parallel positive twist; (c) parallel negative twist; (d) antiparallel positive twist; (e) antiparallel negative twist. \label{f24}}
\end{figure}

\begin{lemma}
By replacing $n$-twist ($n\ge 2$) by $(n+2)$-twist in the Conway symbol of an alternating knot $K$, the signature changes by $-2$ if the replacement is made in a parallel twist with positive crossings, the signature changes by $+2$ if the replacement is made in a parallel twist with negative crossings, and remains unchanged if the replacement is made in an anti-parallel twist.
\end{lemma}

{\bf Proof}: According to the preceding theorem:
\begin{enumerate}
\item by adding a full twist in a parallel positive $n$-twist the writhe changes by $+2$, the number of the white regions $W$ remains unchanged, the number of black regions $B$ increases by $+2$, and the signature changes by $-2$;
\item by adding a full twist in a parallel negative $n$-twist the writhe changes by $-2$, the number of white regions $W$ increases by $2$, the number of black regions $B$ remains unchanged, and the signature increases by 2;
\item  by adding a full twist in an anti-parallel positive $n$-twist the writhe changes by $+2$, the number of white regions $W$ increases by 2, the number of black regions $B$ remains unchanged, and the signature remains unchanged;
\item  by adding a full twist in an anti-parallel negative $n$-twist the writhe changes by $-2$, the number of white regions $W$ remains unchanged, the number of black regions $B$ increases by 2, and the signature remains unchanged.
\end{enumerate}

\begin{figure}[th]
\centerline{\psfig{file=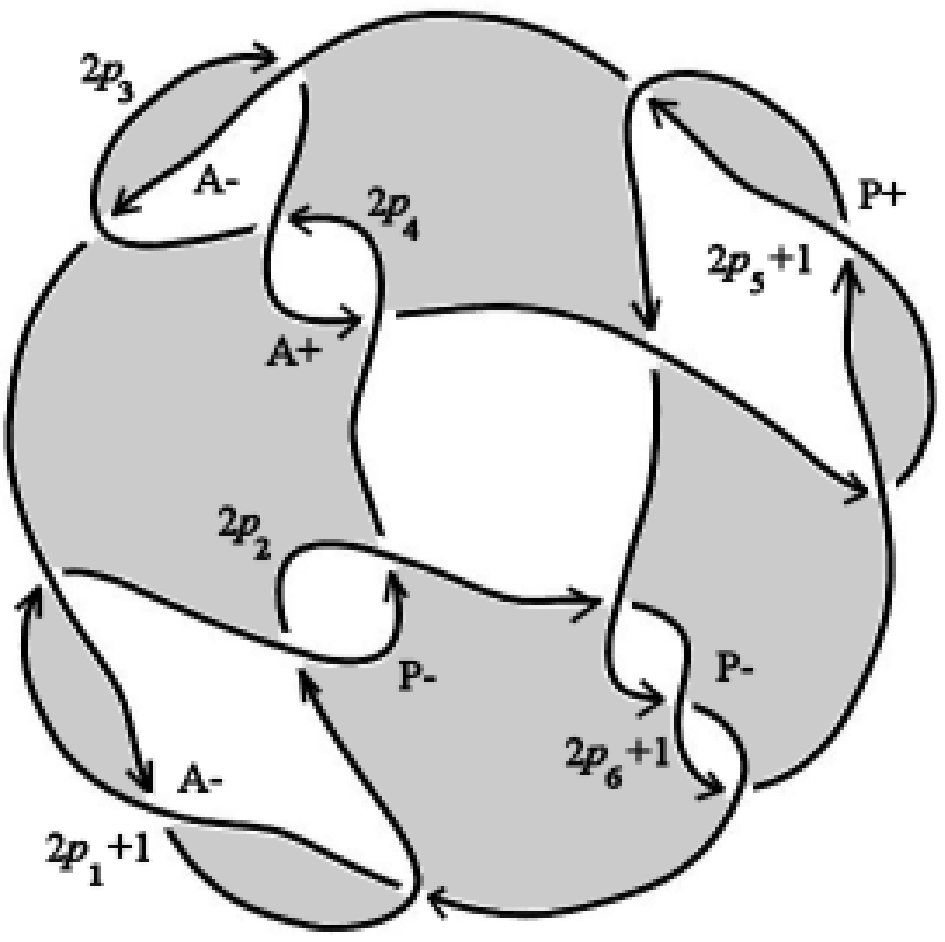,width=1.00in}} \vspace*{8pt}
\caption{Knot family $(2p_1+1)\,(2p_2),(2p_3)\,(2p_4),(2p_5+1)\,1,(2p_6+1)$ beginning with knot $3\,2,2\,2,3\,1,3$. \label{f25}}
\end{figure}

\begin{figure}[th]
\centerline{\psfig{file=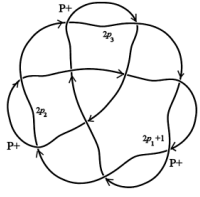,width=1.00in}} \vspace*{8pt}
\caption{Knot family $(2p_1+1):(2p_2):(2p_3)$ beginning with knot $3:2:2$. \label{f26}}
\end{figure}

\begin{figure}[th]
\centerline{\psfig{file=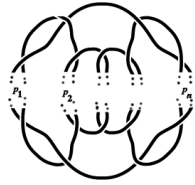,width=1.00in}} \vspace*{8pt}
\caption{Pretzel knot $p_1,p_2,\ldots ,p_n$. \label{f27}}
\end{figure}

\begin{figure}[th]
\centerline{\psfig{file=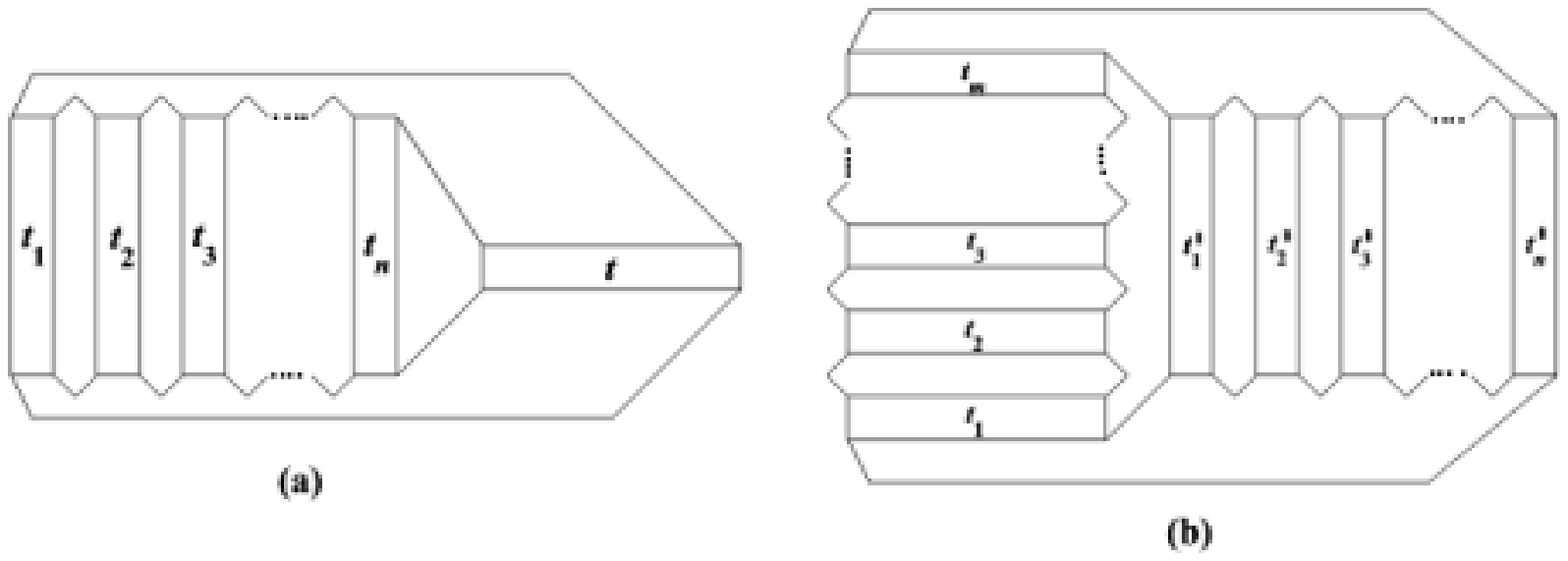,width=4.20in}} \vspace*{8pt}
\caption{(a) Knot $t_1,\ldots t_n+t$; (b) knot $(t_1,\ldots ,t_m)\,(t_1',\ldots ,t_n')$. \label{f28}}
\end{figure}

\begin{theorem}
The signature $\sigma _K$ of an alternating knot $K$ given by its Conway symbol is $$\sigma _K=\sum_P -2[{{n_i} \over 2}]c_i+2c_0,$$ where the sum is taken over all parallel twists $n_i$, $c_i\in \{1,-1\}$ is the sign of crossings belonging to a parallel twist $n_i$, and $2c_0$ is an integer constant which can be computed from the signature of the generating knot.
\end{theorem}

The proof of this theorem follows directly form the preceding Lemma, claiming that only additions of twists in parallel twists in a Conway symbol result in the change of signature, and that by every such addition, signature changes by $-2c_i$. Notice that the condition that we are making twist replacements in the standard Conway symbols, i.e., Conway symbols with the maximal twists, is essential for computation of general formulae for the signature of alternating knot families.

{\bf Example 1}: For the family of Montesinos knots with the Conway symbol of the form $(2p_1+1)\,(2p_2),(2p_3)\,(2p_4),(2p_5+1)\,1,(2p_6+1)$ (Fig. 25), beginning with the generating knot $3\,2,2\,2,3\,1,3$, the parallel twists with negative signs are $2p_2$ and $2p_6+1$, the parallel twist with positive signs is $2p_5+1$, and the remaining twists are anti-parallel. Hence, the signature is $\sigma =-2p_2+2p_5-2p_6+2c_0$. Since the writhe of the generating knot $G=3\,2,2\,2,3\,1,3$ is $w=-4$ and its checkerboard coloring has $W=9$ white and $B=9$ black regions, its signature is $2$. Evaluating the formula $\sigma =2p_2-2p_52p_6+2c_0$ for $\sigma_G= $, $p_2=1$, $p_5=1$, and $p_6=1$, we obtain $c_0=0$. Hence, the general formula for the signature of knots belonging to the family $(2p_1+1)\,(2p_2),(2p_3)\,(2p_4),(2p_5+1)\,1,(2p_6+1)$ is $2p_2-2p_5+2p_6$.

{\bf Example 2}: For the family of polyhedral knots with the Conway symbol of the form $(2p_1+1):(2p_2):(2p_3)$ ($p_1\ge 1$, $p_2\ge 1$, $p_3\ge 1$), beginning with the knot $3:2:2$ (Fig. 26), all twists are parallel twists with positive crossings, and the formula for the signature is $-2p_1-2p_2-2p_3-2$, i.e., $c_0=2$. Constant $c_0$ is computed from the signature of the generating knot $3:2:2$ which is equal to $-4$.

{\bf Example 3}: Let us consider pretzel knots and links (Fig. 27) given by Conway symbol $p_1,\ldots ,p_n$ ($n\ge 3$). We obtain knots if all $p_i$ ($i=1,...,n$) are odd and $n$ is an odd number, or if one twist is even, and all the others are odd. If all twists are odd and $n$ is an odd number, all twists are anti-parallel, and the signature is $\sigma _K=n-1$ for every such knot. If $n=3$, for the pretzel knots of the form $(2p_1+1),(2p_2+1),(2q)$, the twists $2p_1+1$ and $2p_2+1$ are parallel with positive crossings, the twist $2q$ is antiparallel, and the signature is $\sigma _K=2p_1+2p_2$. For $n\ge 4$, for pretzel knots consisting of an even number of odd twists and one even twist, $2p_1+1$, $\ldots $, $2p_{2k}+1$, $2q$, all odd twists are parallel with positive crossings, the even twist $2q$ is anti-parallel, and the signature is $\sigma _K=2p_1+2p_2+\ldots +2p_{2k+1}$. For $n\ge 4$, for pretzel knots consisting of an odd number of odd twists and one even twist, $2p_1+1$, $\ldots $, $2p_{2k+1}+1$, $2q$, all twists are parallel with positive crossings, and the signature is $\sigma _K=2p_1+2p_2+\ldots +2p_{2k+1}+2q$. Hence, for this class of pretzel knots we simply conclude that their unknotting number is given by the formula $u_K=p_1+p_2+\ldots +p_{2k+1}+q.$

{\bf Example 4}: Let us consider knots of the form $t_1,\ldots t_n+t$ ($n\ge 3$), where $t_i$ and $t$ are twists (Fig. 28a).  If the twists of an odd length are denoted by $p$, and twists of an even length by $q$, we have six possible cases:
\begin{enumerate}
\item if the tangle $t_1,\ldots ,t_n$ consists of $2k$ odd twists $p_1,\ldots,p_{2k}$, and the tangle $t$ is an odd twist $p$, the signature is given by the formula $2k+2[{p\over 2}]$
\item if the tangle $t_1,\ldots ,t_n$ consists of $2k+1$ odd twists $p_1,\ldots,p_{2k}$, and the tangle $t$ is an even twist $q$, the signature is given by the formula $2k+q$
\item if the tangle $t_1,\ldots ,t_n$ consists of $2k+1$ odd twists $p_1,\ldots,p_{2k+1}$ and an even twist $q_1$, and the tangle $t$ is an odd twist $p$, the signature is given by the formula $\sum_{i=1}^{2k+1}2[{p_i \over 2}]$
\item if the tangle $t_1,\ldots ,t_n$ consists of $2k+1$ odd twists $p_1,\ldots,p_{2k+1}$ and an even twist $q_1$, and the tangle $t$ is an even twist $q$, the signature is given by the formula $\sum_{i=1}^{2k+1}2[{p_i \over 2}]+q_1$
\item if the tangle $t_1,\ldots ,t_n$ consists of $2k$ odd twists $p_1,\ldots,p_{2k}$ and an even twist $q_1$, and the tangle $t$ is an odd twist $p$, the signature is given by the formula $\sum_{i=1}^{2k}2[{p_i \over 2}]+q_1$
\item if the tangle $t_1,\ldots ,t_n$ consists of $2k$ odd twists $p_1,\ldots,p_{2k}$ and an even twist $q_1$, and the tangle $t$ is an even twist $q$, the signature is given by the formula $\sum_{i=1}^{2k}2[{p_i \over 2}].$
\end{enumerate}

{\bf Example 5}: As a more complex example, we provide general formulae for the signature of knots of the type $(t_1,\ldots ,t_m)\,(t_1',\ldots ,t_n)$ ($m\ge 2$, $n\ge 2$), where twists are denoted by $t_i$ or $t_i'$ (Fig. 28b). If the twists of an odd length are denoted by $p$, and twists of an even length by $q$, we have seven possible cases:

\begin{enumerate}
\item if the first tangle $t_1,\ldots ,t_m$ consists of $2k$ odd twists $p_1,\ldots,p_{2k}$, and the second tangle $t_1',\ldots ,t_n'$ consists of $2r+1$ odd twists $p_1',\ldots,p_{2r+1}'$, the signature is given by the formula
    $$\sum_{i=1}^{2r+1}2[{p_i' \over 2}]+2k$$
\item if the first tangle $t_1,\ldots ,t_m$ consists of $2k$ odd twists $p_1,\ldots,p_{2k}$, and the second tangle $t_1',\ldots ,t_n'$ consists of $2r$ odd twists $p_1',\ldots,p_{2r}'$, the signature is given by the formula
    $$\sum_{i=1}^{2k}2[{p_i \over 2}]-\sum_{i=1}^{2r}2[{p_i' \over 2}]$$
\item if the first tangle $t_1,\ldots ,t_m$ consists of $2k$ odd twists $p_1,\ldots,p_{2k}$ and one even twist $q_1$, and the second tangle $t_1',\ldots ,t_n'$ consists of $2r+1$ odd twists $p_1',\ldots,p_{2r+1}'$, the signature is given by the formula
    $$\sum_{i=1}^{2k}2[{p_i \over 2}]+q_1+2r$$
\item if the first tangle $t_1,\ldots ,t_m$ consists of $2k+1$ odd twists $p_1,\ldots,p_{2k+1}$ and one even twist $q_1$, and the second tangle $t_1',\ldots ,t_n'$ consists of $2r+1$ odd twists $p_1',\ldots,p_{2r+1}'$, the signature is given by the formula
    $$\sum_{i=1}^{2k+1}2[{p_i \over 2}]+2r$$
\item if the first tangle $t_1,\ldots ,t_m$ consists of $2k$ odd twists $p_1,\ldots,p_{2k}$ and one even twist $q_1$, and the second tangle $t_1',\ldots ,t_n'$ consists of $2r$ odd twists $p_1',\ldots,p_{2r}'$ and one even twist $q_1'$, the signature is given by the formula
    $$\sum_{i=1}^{2k}2[{p_i \over 2}]-\sum_{i=1}^{2r}2[{p_i' \over 2}]$$
\item if the first tangle $t_1,\ldots ,t_m$ consists of $2k$ odd twists $p_1,\ldots,p_{2k}$ and one even twist $q_1$, and the second tangle $t_1',\ldots ,t_n'$ consists of $2r+1$ odd twists $p_1',\ldots,p_{2r+1}'$ and one even twist $q_1'$, the signature is given by the formula
    $$\sum_{i=1}^{2k}2[{p_i \over 2}]-\sum_{i=1}^{2r+1}2[{p_i' \over 2}]-q_1'$$
\item if the first tangle $t_1,\ldots ,t_m$ consists of $2k+1$ odd twists $p_1,\ldots,p_{2k+1}$ and one even twist $q_1$, and the second tangle $t_1',\ldots ,t_n'$ consists of $2r+1$ odd twists $p_1',\ldots,p_{2r+1}'$ and one even twist $q_1'$, the signature is given by the formula
    $$\sum_{i=1}^{2k+1}2[{p_i \over 2}]-\sum_{i=1}^{2r+1}2[{p_i' \over 2}]+q_1-q_1'.$$
\end{enumerate}

\section{Unknotting numbers of knot families}

K.~Murasugi \cite{9} proved the lower bound for the unknotting number of knots, $u(K)\geq  \frac{|\sigma _K|}{2}$. Using this criterion, for many (sub)families of knots  we can confirm that their $BJ$-unknotting numbers, i.e., unknotting numbers computed according to Bernhard-Jablan Conjecture \cite{5}
represent the actual unknotting numbers of these (sub)families. In the following table is given the list of (sub)families with this property obtained from knots with at most $n=8$ crossings, where in the first column is given the first knot belonging to the family, in the second its Conway symbol, in the third the general Conway symbol, in the fourth the general formula for the signature, in the fifth the unknotting number confirmed by the signature, and in the sixth the conditions for this unknotting number\footnote{Conditions for unknotting numbers are determined from the experimental results obtained for knots up to $n=20$ crossings.}.

\medskip

\scriptsize

\noindent \begin{tabular}{|@{$\,$}c@{$\,$}|@{$\,$}c@{$\,$}|@{$\,$}c@{$\,$}|@{$\,$}c@{$\,$}|@{$\,$}c@{$\,$}|@{$\,$}c@{$\,$}|} \hline

$K$ & $Con$ & $Fam$ & $\sigma $ & $u$ & $Cond$   \\ \hline
$3_1$ & $3$ & $(2p_1+1)$ & $2p_1$ & $p_1$ & $$   \\ \hline
$4_1$ & $2\,2$ & $(2p_1)\,(2p_2)$ & $0$ & $$ & $$   \\ \hline
$5_2$ & $3\,2$ & $(2p_1+1)\,(2p_2)$ & $2p_2$ & $p_2$ & $$   \\ \hline
$6_2$ & $3\,1\,2$ & $(2p_1+1)\,1\,(2p_2)$ & $2p_1$ & $p_1$ & $p_1\ge p_2$   \\ \hline
$6_3$ & $2\,1\,1\,2$ & $(2p_1)\,1\,1\,(2p_2)$ & $2p_1-2p_2$ & $|p_1-p_2|$ & $p_1\neq p_2$   \\ \hline
$7_4$ & $3\,1\,3$ & $(2p_1+1)\,1\,(2p_2+1)$ & $2$ & $$ & $$   \\ \hline
$7_5$ & $3\,2\,2$ & $(2p_1+1)\,(2p_2)\,(2p_3)$ & $2p_1+2p_3$ & $p_1+p_3$ & $$   \\ \hline
$7_6$ & $2\,2\,1\,2$ & $(2p_1)\,(2p_2)\,1\,(2p_3)$ & $2p_3$ & $p_3$ & $p_2\le p_3$   \\ \hline
$7_7$ & $2\,1\,1\,1\,2$ & $(2p_1)\,1\,1\,1\,(2p_2)$ & $0$ & $$ & $$   \\ \hline
$8_5$ & $3,3,2$ & $(2p_1+1),(2p_2+1),(2p_3)$ & $2p_1+2p_2$ & $p_1+p_2$ & $p_1\ge p_3$ or $p_2\ge p_3$   \\ \hline
$8_6$ & $3\,3\,2$ & $(2p_1+1)\,(2p_2+1)\,(2p_3)$ & $2p_1$ & $$ & $$   \\ \hline
$8_8$ & $2\,3\,1\,2$ & $(2p_1)\,(2p_2+1)\,1\,(2p_3)$ & $2p_1-2p_3$ & $p_3-p_1$ & $p_3-p_1>p_2$   \\ \hline
$8_9$ & $3\,1\,1\,3$ & $(2p_1+1)\,1\,1\,(2p_2+1)$ & $2p_1-2p_2$ & $|p_1-p_2|$ & $p_1\neq p_2$   \\ \hline
$8_{10}$ & $3,2\,1,2$ & $\,\,(2p_1+1),(2p_2)\,1,(2p_3)\,\,$ & $\,\,2p_1-2p_2+2p_3\,\,$ & $p_1-p_2+p_3$ & $p_3>p_2$   \\ \hline
$8_{11}$ & $3\,2\,1\,2$ & $(2p_1+1)\,2p_2)\,1\,(2p_3)$ & $2p_2$ & $p_2$ & $p_2\ge p_3$   \\ \hline
$8_{12}$ & $2\,2\,2\,2$ & $(2p_1)\,(2p_2)\,(2p_3)\,(2p_4)$ & $0$ & $$ & $$   \\ \hline
$8_{13}$ &  $3\,1\,1\,1\,2$ & $(2p_1+1)\,1\,1\,1\,(2p_2)$ & $2p_2-2$ & $p_2-1$ & $p_2-1>p_1$   \\ \hline
$8_{14}$ & $2\,2\,1\,1\,2$ & $(2p_1)\,(2p_2)\,1\,1\,(2p_3)$ & $2p_1$ & $p_1$ & $p_2\le p_3$   \\ \hline
$8_{15}$ & $2\,1,2\,1,2$ & $(2p_1)\,1,(2p_2)\,1,(2p_3)$ & $2p_1+2p_2$ & $p_1+p_2$ & $$   \\ \hline
$8_{16}$ & $.2.2\,0$ & $.(2p_1).(2p_2)\,0$ & $2p_1+2p_2-2$ & $p_1+p_2-1$ & $$   \\ \hline
$8_{17}$ & $.2.2$ & $.(2p_1).(2p_2)$ & $2p_1-2p_2$ & $|p_1-p_2|$ & $p_1=1,p_2>1$   \\
$$ & $$ & $$ & $$ & $$ & or $p_2=1,p_1>1$   \\ \hline

\end{tabular}

\scriptsize

\noindent \begin{tabular}{|@{$\,$}c@{$\,$}|@{$\,$}c@{$\,$}|@{$\,$}c@{$\,$}|@{$\,$}c@{$\,$}|@{$\,$}c@{$\,$}|@{$\,$}c@{$\,$}|} \hline
$K$ & $Con$ & $Fam$ & $\sigma $ & $u$ & $Cond$   \\ \hline
$9_{10}$ & $3\,3\,3$ & $(2p_1+1)\,(2p_2+1)\,(2p_3+1)$ & $2p_2+2$ & $$ & $$   \\ \hline
$9_{13}$ & $3\,2\,1\,3$ & $(2p_1)\,(2p_2)\,1\,(2p_3+1)$ & $2p_1+2$ & $$ & $$   \\ \hline
$9_{15}$ & $2\,3\,2\,2$ & $(2p_1)\,(2p_2+1)\,(2p_3)\,(2p_4)$ & $2p_1$ & $$ & $$   \\ \hline
$9_{16}$ & $3,3,2+$ & $(2p_1+1),(2p_2+1),(2p_3)+$ & $2p_1+2p_2+2p_3$ & $p_1+p_2+p_3$ & $$   \\ \hline
$9_{17}$ & $2\,1\,3\,1\,2$ & $(2p_1)\,1\,(2p_2+1)\,1\,(2p_3)$ & $-2p_2$ & $p_2$ & $p_1+p_3\le p_2$   \\ \hline
$9_{18}$ & $3\,2\,2\,2$ & $(2p_1+1)\,(2p_2)\,(2p_3)\,(2p_4)$ & $2p_2+2p_4$ & $p_2+p_4$ & $$   \\ \hline
$9_{19}$ & $2\,3\,1\,1\,2$ & $(2p_1)\,(2p_2+1)\,1\,1\,(2p_3)$ & $0$ & $$ & $$   \\ \hline
$9_{20}$ & $3\,1\,2\,1\,2$ & $(2p_1+1)\,1\,(2p_2)\,1\,(2p_3)$ & $2p_1+2p_3$ & $p_1+p_3$ & $p_1+p_3\ge p_2$   \\ \hline
$9_{21}$ & $3\,1\,1\,2\,2$ & $(2p_1+1)\,1\,1\,(2p_2)\,(2p_3)$ & $2$ & $$ & $$   \\ \hline
$9_{22}$ & $2\,1\,1,3,2$ & $(2p_1)\,1\,1,(2p_2+1),(2p_3)$ & $-2p_2$ & $p_2$ & $p_1+p_3-1\le p_2$   \\ \hline
$9_{23}$ & $2\,2\,1\,2\,2$ & $(2p_1)\,(2p_2)\,1\,(2p_3)\,(2p_4)$ & $2p_1+2p_4$ & $p_1+p_4$ & $$   \\ \hline
$9_{24}$ & $2\,1,3,2+$ & $(2p_1)\,1,(2p_2+1),(2p_3)+$ & $2p_1-2p_2$ & $$ & $$   \\ \hline
$9_{25}$ & $2\,2,2\,1,2$ & $(2p_1)\,(2p_2),(2p_3)\,1,(2p_4)$ & $-2p_3$ & $$ & $$   \\ \hline
$9_{26}$ & $3\,1\,1\,1\,1\,2$ & $(2p_1+1)\,1\,1\,1\,1\,(2p_2)$ & $2p_1$ & $p_1$ & $$   \\ \hline
$9_{27}$ & $2\,1\,2\,1\,1\,2$ & $(2p_1)\,1\,(2p_2)\,1\,1\,(2p_3)$ & $2p_3-2p_2$ & $p_3-p_2$ & $p_2<p_3$   \\ \hline
$9_{28}$ & $2\,1,2\,1,2+$ & $(2p_1)\,1,(2p_2)\,1,(2p_3)+$ & $2p_1+2p_2-2p_3$ & $p_1+p_2-p_3$ & $p_3\le p_1$ or $p_3\le p_2$   \\ \hline
$9_{29}$ & $.2.2\,0.2$ & $.(2p_1).(2p_2)\,0.(2p_3)$ & $-2p_2$ & $$ & $$   \\ \hline
$9_{30}$ & $2\,1\,1,2\,1,2$ & $(2p_1)\,1\,1,(2p_2)\,1,(2p_3)$ & $2p_2-2p_3$ & $p_3-p_2$ & $p_1+p_2\le p_3$   \\ \hline
$9_{31}$ & $2\,1\,1\,1\,1\,1\,2$ & $(2p_1)\,1\,1\,1\,1\,1\,(2p_2)$ & $2p_1-2p_2-2$ & $p_1+p_2-1$ & $p_1+p_2>2$   \\ \hline
$9_{32}$ & $.2\,1.2\,0$ & $.(2p_1)\,1.(2p_2)\,0$ & $2p_2$ & $$ & $$   \\ \hline
$9_{33}$ & $.2\,1.2$ & $.(2p_1)\,1.(2p_2)$ & $-2p_2+2$ & $p_2-1$ & $p_2-p_1\ge 2$   \\ \hline
$9_{34}$ & $8^*2\,0$ & $8^8(2p_1)\,0$ & $0$ & $$ & $$   \\ \hline
$9_{35}$ & $3,3,3$ & $(2p_1+1),(2p_2+1),(2p_3+1)$ & $2$ & $$ & $$   \\ \hline
$9_{36}$ & $2\,2,3,2$ & $(2p_1)\,(2p_2),(2p_3)\,1,(2p_4)$ & $2p_3+2p_4$ & $p_3+p_4$ & $p_2\le p_4$   \\ \hline
$9_{37}$ & $2\,1,2\,1,3$ & $(2p_1)\,1,(2p_2)\,1,(2p_3+1)$ & $0$ & $$ & $$   \\ \hline
$9_{38}$ & $.2.2.2$ & $.(2p_1).(2p_2).(2p_3)$ & $2p_2+2$ & $$ & $$   \\ \hline
$9_{39}$ & $2:2:2\,0$ & $(2p_1):(2p_2):(2p_3)\,0$ & $2$ & $$ & $$   \\ \hline
$9_{41}$ & $2\,0:2\,0:2\,0$ & $(2p_1)\,0:(2p_2)\,0:(2p_3)\,0$ & $0$ & $$ & $$   \\ \hline
\end{tabular}

\bigskip

\begin{landscape}

\scriptsize

\noindent \begin{tabular}{|@{$\,$}c@{$\,$}|@{$\,$}c@{$\,$}|@{$\,$}c@{$\,$}|@{$\,$}c@{$\,$}|@{$\,$}c@{$\,$}|@{$\,$}c@{$\,$}|} \hline
$K$ & $Con$ & $Fam$ & $\sigma $ & $u$ & $Cond$   \\ \hline
$10_{22}$ & $3\,3\,1\,3$ & $(2p_1+1)\,(2p_2+1)\,1\,(2p_1+1)$ & $2p_1-2p_3$ & $p_3-p_1$ & $p_3-p_1>p_2$    \\ \hline
$10_{23}$ & $3\,3\,1\,1\,2$ & $(2p_1+1)\,(2p_2+1)\,1\,1,(2p_3)$ & $2p_2-2p_3+2$ & $p_2-p_3+1$  & $p_3\le p_2$  \\ \hline
$10_{24}$ & $3\,2\,3\,2$ & $(2p_1+1)\,(2p_2)\,(2p_3+1)\,(2p_4)$ & $2p_2$ & $$  & $$  \\ \hline
$10_{25}$ & $3\,2\,2\,1\,2$ & $(2p_1+1)\,(2p_2)\,(2p_3)\,1\,(2p_4)$ & $2p_1+2p_3$ & $p_1+p_3$  & $p_4\le p_3$  \\ \hline
$10_{26}$ & $3\,2\,1\,1\,3$ & $(2p_1+1)\,(2p_2)\,1\,1\,(2p_3+1)$ & $2p_2-2p_3$ & $|p_2-p_3|$  & $p_3-p_2>p_1$ or $p_2>p_3$  \\ \hline
$10_{27}$ & $3\,2\,1\,1\,1\,2$ & $(2p_1+1)\,(2p_2)\,1\,1\,1(2p_3)$ & $2p_1-2p_3+2$ & $|p_1-p_3+1|$  & $p_3\le p_1, p_2=2$  \\
$$ & $$ & $$ & $$ & $$  & or $p_3-p_1>p_2+1$  \\ \hline
$10_{28}$ & $3\,1\,3\,1\,2$ & $(2p_1+1)\,1\,(2p_2+1)\,1\,(2p_3)$ & $-2p_3+2$ & $p_3-1$  & $p_1+p_2+1<p_3$  \\ \hline
$10_{29}$ & $3\,1\,2\,2\,2$ & $(2p_1+1)\,1\,(2p_2)\,(2p_3)\,(2p_4)$ & $2p_1$ & $$  & $$  \\ \hline
$10_{30}$ & $3\,1\,2\,1\,\,1\,2$ & $(2p_1+1)\,1\,(2p_2)\,1\,1\,(2p_3)$ & $2$ & $$  & $$  \\ \hline
$10_{31}$ & $3\,1\,1\,3\,2$ & $(2p_1+1)\,1\,1\,(2p_2+1)\,(2p_3)$ & $-2p_3+2$ & $$  & $$  \\ \hline
$10_{32}$ & $3\,1\,1\,1\,2\,2$ & $(2p_1+1)\,1\,1\,1\,(2p_2)\,(2p_3)$ & $2p_1-2p_3$ & $|p_1-p_3|$  & $p_1>p_2+p_3$  \\
$$ & $$ & $$ & $$ & $$  & or $p_2=2, p_3>p_1$  \\ \hline
$10_{33}$ & $3\,1\,1\,1\,1\,3$ & $(2p_1+1)\,1\,1\,1\,1\,(2p_2+1)$ & $0$ & $$  & $$  \\ \hline
$10_{37}$ & $2\,3\,3\,2$ & $(2p_1)\,(2p_2+1)\,(2p_3+1)\,(2p_4)$ & $2p_1-2p_4$ & $$  & $$  \\ \hline
$10_{38}$ & $2\,3\,1\,2\,2$ & $(2p_1)\,(2p_2+1)\,1\,(2p_3)\,(2p_4)$ & $-2p_4$ & $$  & $$  \\ \hline
$10_{39}$ & $2\,2\,3\,1\,2$ & $(2p_1)\,(2p_2)\,(2p_3+1)\,1\,(2p_4)$ & $2p_1+2p_3$ & $p_1+p_3$  & $p_1+p_2+p_3\le p_4$  \\ \hline
$10_{40}$ & $2\,2\,2\,1\,1\,2$ & $(2p_1)\,(2p_2)\,(2p_3)\,1\,1\,(2p_4)$ & $2p_1+2p_3-2p_4$ & $|p_1+p_3-p_4|$  & $p_4<p_3$  \\
$$ & $$ & $$ & $$ & $$  & or $p_4>p_1+p_2+p_3$  \\ \hline
$10_{41}$ & $2\,2\,1\,2\,1\,2$ & $(2p_1)\,(2p_2)\,1\,(2p_3)\,1\,(2p_4)$ & $2p_3$ & $p_3$  & $p_2+p_4\le p_3$  \\ \hline
$10_{42}$ & $2\,2\,1\,1\,1\,1\,2$ & $(2p_1)\,(2p_2)\,1\,1\,1\,1\,(2p_3)$ & $-2p_3+2$ & $$  & $$  \\ \hline
$10_{43}$ & $2\,1\,2\,2\,1\,2$ & $(2p_1)\,1\,(2p_2)\,(2p_3)\,1\,(2p_4)$ & $2p_1-2p_4$ & $$  & $$  \\ \hline
$10_{44}$ & $2\,1\,2\,1\,1\,1\,2$ & $(2p_1)\,1\,(2p_2)\,1\,1\,1\,(2p_3)$ & $2p_1$ & $p_1$  & $p_2+p_3\le p_1+1$  \\ \hline
$10_{45}$ & $2\,1\,1\,1\,1\,1\,1\,2$ & $(2p_1)\,1\,1\,1\,1\,1\,1\,(2p_2)$ & $0$ & $$  & $$  \\ \hline
$10_{50}$ & $3\,2,3,2$ & $(2p_1+1)\,(2p_2),(2p_3+1),(2p_4)$ & $2p_2+2p_3$ & $p_2+p_3$  & $p_4\le p_2$  \\ \hline
$10_{51}$ & $3\,2,2\,1,2$ & $(2p_1+1)\,(2p_2),(2p_3)\,1,(2p_4)$ & $2p_2+2p_4-2p_3$ & $$  & $$  \\ \hline
$10_{52}$ & $3\,1\,1,3,2$ & $(2p_1+1)\,1\,1,(2p_2+1),(2p_3)$ & $-2p_2-2p_3+2$ & $p_2+p_3-1$  & $p_3-p_1\ge 2$  \\ \hline
$10_{53}$ & $3\,1\,1,2\,1,2$ & $(2p_1+1)\,1\,1,(2p_2)\,1,(2p_3)$ & $2p_2+2$ & $$  & $$  \\ \hline
$10_{54}$ & $2\,3,3,2$ & $(2p_1)\,(2p_2+1),(2p_3+1),(2p_4)$ & $2p_1-2p_3-2p_4$ & $p_3+p_4-p_1$  & $p_4>p_1+p_2$  \\ \hline
$10_{55}$ & $2\,3,2\,1,2$ & $(2p_1)\,(2p_2+1),(2p_3)\,1,(2p_4)$ & $2p_1+2p_3$ & $p_1+p_3$  & $$  \\ \hline
$10_{56}$ & $2\,2\,1,3,2$ & $(2p_1)\,(2p_2)\,1,(2p_3+1),(2p_4)$ & $2p_1+2p_3$ & $p_1+p_3$  & $p_4\le p_1+p_2$  \\ \hline
$10_{57}$ & $2\,2\,1,2\,1,2$ & $(2p_1)\,(2p_2)\,1,(2p_3)\,1,(2p_4)$ & $2p_1+2p_4-2p_3$ & $p_1+p_4-p_3$  & $p_3<p_4$  \\ \hline
$10_{58}$ & $2\,1\,1,2\,1\,1,2$ & $(2p_1)\,1\,1,(2p_2)\,1\,1,(2p_3)$ & $0$ & $$  & $$  \\ \hline
$10_{59}$ & $2\,2,2\,1\,1,2$ & $(2p_1)\,(2p_2),(2p_3)\,1\,1,(2p_4)$ & $2p_4$ & $p_4$  & $p_2+p_3-1\le p_4$  \\ \hline
$10_{60}$ & $2\,1\,1,2\,1\,1,2$ & $(2p_1)\,1\,1,(2p_2)\,1\,1,(2p_3)$ & $0$ & $$  & $$  \\ \hline
$10_{64}$ & $3\,1,3,3$ & $(2p_1+1)\,1,(2p_2+1),(2p_3+1)$ & $2p_1-2p_2-2p_3$ & $p_2+p_3-p_1$  & max$(p_2,p_3)>p_1$  \\ \hline

\end{tabular}
\end{landscape}

\bigskip

\begin{landscape}

\scriptsize

\noindent \begin{tabular}{|c|c|c|c|c|c|} \hline
$K$ & $Con$ & $Fam$ & $\sigma $ & $u$ & $Cond$   \\ \hline
$10_{65}$ & $3\,1,3,2\,1$ & $(2p_1+1)\,1,(2p_2+1),(2p_3)\,1$ & $2p_2-2p_3+2$ & $$  & $$  \\ \hline
$10_{66}$ & $3\,1,2\,1,2\,1$ & $(2p_1+1)\,1,(2p_2)\,1,(2p_3)\,1$ & $2p_1+2p_2+2p_3$ & $p_1+p_2+p_3$  & $$  \\ \hline
$10_{68}$ & $2\,1\,1,3,3$ & $(2p_1)\,1\,1,(2p_2+1),(2p_3+1)$ & $2p_1-2$ & $$  & $$  \\ \hline
$10_{67}$ & $2\,2,3,2\,1$ & $(2p_1)\,(2p_2),(2p_3+1),(2p_4)\,1$ & $2p_1$ & $$  & $$  \\ \hline
$10_{69}$ & $2\,1\,1,2\,1,2\,1$ & $(2p_1)\,1\,1,(2p_2)\,1,(2p_3)\,1$ & $2p_1$ & $$  & $$  \\ \hline
$10_{70}$ & $2\,2,3,2+$ & $(2p_1)\,(2p_2),(2p_3+1),(2p_4)+$ & $2p_3$ & $$  & $$  \\ \hline
$10_{71}$ & $2\,2,2\,1,2+$ & $(2p_1)\,(2p_2),(2p_3)\,1,(2p_4)+$ & $2p_4-2p_3$ & $$  & $$  \\ \hline
$10_{72}$ & $2\,1\,1,3,2+$ & $(2p_1)\,1\,1,(2p_2+1),(2p_3)+$ & $-2p_2-2p_3$ & $p_2+p_3$  & $p_1-1\le p_2+p_3$  \\ \hline
$10_{73}$ & $2\,1\,1,2\,1,2+$ & $(2p_1)\,1\,1,(2p_2)\,1,(2p_3)+$ & $2p_2$ & $p_2$  & $p_3=1$  \\ \hline
$10_{74}$ & $3,3,2\,1+$ & $(2p_1+1),(2p_2+1),(2p_3)\,1+$ & $2$ & $$  & $$  \\ \hline
$10_{75}$ & $2\,1,2\,1,2\,1+$ & $(2p_1)\,1,(2p_2)\,1,(2p_3)\,1+$ & $0$ & $$  & $$  \\ \hline
$10_{76}$ & $3,3,2+2$ & $(2p_1+1),(2p_2+1),(2p_3)+(2p_4)$ & $2p_1+2p_2$ & $$  & $$  \\ \hline
$10_{77}$ & $3,2\,1,2+2$ & $(2p_1+1),(2p_2)\,1,(2p_3)+(2p_4)$ & $2p_1+2p_3-2p_2$ & $$  & $$  \\ \hline
$10_{78}$ & $2\,1,2\,1,2+2$ & $(2p_1)\,1,(2p_2)\,1,(2p_3)+(2p_4)$ & $2p_1+2p_2$ & $p_1+p_2$  & $p_4\le p_1+p_2$  \\ \hline
$10_{79}$ & $(3,2)\,(3,2)$ & $((2p_1+1),(2p_2))\,((2p_3+1),(2p_4))$ & $2p_1+2p_2-2p_3-2p_4$ & $$  & $$  \\ \hline
$10_{80}$ & $(3,2)\,(2\,1,2)$ & $((2p_1+1),(2p_2))\,((2p_3)\,1,(2p_4))$ & $2p_1+2p_2+2p_3$ & $p_1+p_2+p_3$  & $$  \\ \hline
$10_{81}$ & $(2\,1,2)\,(2\,1,2)$ & $((2p_1)\,1,(2p_2))\,((2p_3)\,1,(2p_4))$ & $2p_1-2p_3$ & $$  & $$  \\ \hline
$10_{83}$ & $.3\,1.2$ & $.(2p_1+1)\,1.(2p_2)$ & $-2p_2+2$ & $p_2-1$  & $p_2>p_1+1$  \\ \hline
$10_{84}$ & $.2\,2.2$ & $.(2p_1)\,(2p_2).(2p_3)$ & $2p_1+2p_3$ & $p_1+p_3$  & $p_2=1$ or $p_3\ge 2$  \\ \hline
$10_{86}$ & $.3\,1.2\,0$ & $.(2p_1+1)\,1.(2p_2)\,0$ & $2p_2$ & $$  & $$  \\ \hline
$10_{87}$ & $.2\,2.2\,0$ & $.(2p_1)\,(2p_2).(2p_3)\,0$ & $2p_1-2p_3$ & $$  & $$  \\ \hline
$10_{88}$ & $.2\,1.2\,1$ & $.(2p_1)\,1.(2p_2)\,1$ & $0$ & $$  & $$  \\ \hline
$10_{89}$ & $.2\,1.2\,1\,0$ & $.(2p_1)\,1.(2p_2)\,1\,0$ & $2$ & $$  & $$  \\ \hline
$10_{90}$ & $.3.2.2$ & $.(2p_1+1).(2p_2).(2p_3)$ & $2p_2-2p_3$ & $p_3-p_2$  & $p_3>p_1+p_2$  \\ \hline
$10_{91}$ & $.3.2.2\,0$ & $.(2p_1+1).(2p_2).(2p_3)\,0$ & $2p_1-2p_2-2p_3+2$ & $p_2+p_3-p_1-1$  & $p_1\le p_3, p_2>1$  \\ \hline
$10_{92}$ & $.2\,1.2.2\,0$ & $.(2p_1)\,1.(2p_2).(2p_3)\,0$ & $2p_1+2p_2$ & $p_1+p_2$  & $p_3-1\le p_1$  \\ \hline
$10_{93}$ & $.3.2\,0.2$ & $.(2p_1+1).(2p_2)\,0.(2p_3)$ & $-2p_2-2p_3$ & $p_2+p_3$  & $p_1<p_2$ or $p_3>p_1+1$  \\ \hline
$10_{94}$ & $.3\,0.2.2$ & $.(2p_1+1)\,0.(2p_2).(2p_3)$ & $2p_1+2p_2-2p_3$ & $p_1+p_2-p_3$  & $p_1-p_3\ge 1$ or $p_2-p_3\ge 1$ \\ \hline
$10_{95}$ & $.2\,1\,0.2.2$ & $(2p_1)\,1\,0.(2p_2).(2p_3)$ & $2p_1-2p_2-2$ & $p_2-p_1+1$ &  $p_1=p_3=1$  \\ \hline
$10_{96}$ & $.2.2\,1.2$ & $.(2p_1).(2p_2)\,1.(2p_3)$ & $0$ & $$  & $$  \\ \hline
$10_{97}$ & $.2.2\,1\,0.2$ & $.(2p_1).(2p_2)\,1\,0.(2p_3)$ & $2$ & $$  & $$  \\ \hline
$10_{98}$ & $.2.2.2.2\,0$ & $.(2p_1).(2p_2).(2p_3).(2p_4)\,0$ & $2p_1+2p_3$ & $p_1+p_3$  & $p_4\le p_1$ or $p_4\le p_3$  \\ \hline
$10_{99}$ & $.2.2.2\,0.2\,0$ & $.(2p_1).(2p_2).(2p_3)\,0.(2p_4)\,0$ & $2p_1+2p_4-2p_2-2p_3$ & $$  & $$  \\ \hline
$10_{100}$ & $3:2:2$ & $(2p_1+1):(2p_2):(2p_3)$ & $2p_1+2p_2+2p_3-2$ & $p_1+p_2+p_3-1$  & $p_2>1$ or $p_3>1$  \\ \hline
$10_{101}$ & $2\,1:2:2$ & $(2p_1)\,1:(2p_2):(2p_3)$ & $2p_1+2$ & $$  & $$  \\ \hline
$10_{102}$ & $3:2:2\,0$ & $(2p_1+1):(2p_2):(2p-3)\,0$ & $2p_3-2p_2$ & $p_2-p_3$  & $p_3=1,p_2-p_1>1$  \\ \hline
$10_{103}$ & $3\,0:2:2$ & $(2p_1+1)\,0:(2p_2):(2p_3)$ & $2p_2+2p_3-2$ & $$  & $$  \\ \hline
$10_{104}$ & $3:2\,0:2\,0$ & $(2p_1+1):(2p_2)\,0:(2p_3)\,0$ & $2p_1-2p_2-2p_3+2$ & $p_1-p_2-p_3+1$  & $p_1>p_2, p_3=1$  \\
$$ & $$ & $$ & $$ & $$  & or $p_1>p_3, p_2=1$  \\ \hline

\end{tabular}

\bigskip

\noindent \begin{tabular}{|c|c|c|c|c|c|} \hline
$K$ & $Con$ & $Fam$ & $\sigma $ & $u$ & $Cond$   \\ \hline
$10_{105}$ & $2\,1:2\,0:2\,0$ & $(2p_1)\,1:(2p_2)\,0:(2p_3)\,0$ & $2p_1$ & $p_1$  & $p_1>p_2+p_3$  \\ \hline
$10_{106}$ & $3\,0:2:2\,0$ & $(2p_1+1)\,0:(2p_2):(2p_3)\,0$ & $2p_1+2p_3-2p_2$ & $p_1+p_3-p_2$  & $p_2+1\le p_1$  \\ \hline
$10_{107}$ & $2\,1\,0:2:2\,0$ & $(2p_1)\,1\,0:(2p_2):(2p_3)\,0$ & $2p_1-2$ & $$  & $$  \\ \hline
$10_{108}$ & $3\,0:2\,0:2\,0$ & $(2p_1+1)\,0:(2p_2)\,0:(2p_3)\,0$ & $-2p_2-2p_3+2$ & $p_2+p_3-1$  & $p_2>p_1+2$ or $p_3>p_1+2$  \\ \hline
$10_{109}$ & $2.2.2.2$ & $(2p_1).(2p_2).(2p_3).(2p_4)$ & $2p_1+2p_3-2p_2-2p_4$ & $p_1+p_3-p_2-p_4$  & $p_1\ge p_2+p_4$  \\ \hline
$10_{110}$ & $2.2.2.2\,0$ & $(2p_1).(2p_2).(2p_3).(2p_4)\,0$ & $-2p_4$ & $p_4$  & $p_4\ge p_1+p_3$  \\ \hline
$10_{111}$ & $2.2.2\,0.2$ & $(2p_1).(2p_2).(2p_3)\,0.(2p_4)$ & $2p_2+2p_3$ & $p_2+p_3$  & $p_2+p_3\ge p_4\ge p_1$  \\ \hline
$10_{112}$ & $8^*3$ & $8^*(2p_1+1)$ & $2p_1$ & $p_1$  & $p_1\ge 2$ \\ \hline
$10_{113}$ & $8^*2\,1$ & $8^*(2p_1)\,1$ & $2p_1$ & $p_1$  & $p_1\ge 2$  \\ \hline
$10_{114}$ & $8^*3\,0$ & $8^*(2p_1+1)\,0$ & $0$ & $$  & $$  \\ \hline
$10_{115}$ & $8^*2\,0.2\,0$ & $8^*(2p_1)\,0.(2p_2)\,0$ & $0$ & $$  & $$  \\ \hline
$10_{116}$ & $8^*2:2$ & $8^*(2p_1):(2p_2)$ & $2p_1+2p_2-2$ & $p_1+p_2-1$  & $p_1\ge 2$ or $p_2\ge 2$  \\ \hline
$10_{117}$ & $8^*2:2\,0$ & $8^*(2p_1):(2p_2)\,0$ & $2p_2$ & $$  & $$  \\ \hline
$10_{118}$ & $8^*2:.2$ & $8^*(2p_1):.(2p_2)\,0$ & $2p_1-2p_2$ & $|p_1-p_2|$  & $p_1\ge 2,p_2=1$  \\
$$ & $$ & $$ & $$ & $$  & $p_2\ge 2,p_1=1$ or $|p_1-p_2|\ge 2$ \\ \hline
$10_{119}$ & $8^*2:.2\,0$ & $8^*(2p_1):.(2p_2)\,0$ & $-2p_2+2$ & $p_2-1$  & $p_2-p_1\ge 2$  \\ \hline
$10_{120}$ & $8^*2\,0::2\,0$ & $8^*(2p_1)\,0::(2p_2)\,0$ & $4$ & $$  & $$  \\ \hline
$10_{121}$ & $9^*2\,0$ & $9^8(2p_1)\,0$ & $2$ & $$  & $$  \\ \hline
$10_{122}$ & $9^*.2\,0$ & $9^*.(2p_1)\,0$ & $0$ & $$  & $$  \\ \hline
\end{tabular}

\bigskip

\normalsize

In the following table we provide the same results for link families obtained from generating links with at most $n=9$ crossings.

\medskip

\scriptsize

\noindent \begin{tabular}{|c|c|c|c|c|c|} \hline
$K$ & $Con$ & $Fam$ & $\sigma $ & $u$ & $Cond$   \\ \hline
$2_1^2$ & $2$ & $(2p_1)$ & $-2p_1+1$ & $p_1$ & $$   \\ \hline
$5_1^2$ & $2\,1\,2$ & $(2p_1)\,1\,(2p_2)$ & $-2p_1+1$ & $p_1$ & $p_1>p_2$   \\ \hline
$6_2^2$ & $3\,3$ & $(2p_1+1)\,(2p_2+1)$ & $-2p_1-1$ & $$ & $$   \\ \hline
$6_3^2$ & $2\,2\,2$ & $(2p_1)\,(2p_2)\,(2p_3)$ & $-2p_1-2p_3+1$ & $p_1+p_3$ & $$   \\ \hline
$7_2^2$ & $3\,1\,1\,2$ & $(2p_1+1)\,1\,1\,(2p_2)$ & $-2p_1+2p_2$ & $|p_1-p_2|$ & $$   \\ \hline
$7_3^2$ & $2\,3\,2$ & $(2p_1)\,(2p_2+1)\,(2p_3)$ & $-2p_1+1$ & $$ & $$   \\ \hline
$7_4^2$ & $3,2,2$ & $(2p_1+1),(2p_2),(2p_3)$ & $-2p_1-2p_3+1$ & $p_1+p_3$ & $p_2=1$   \\ \hline
$7_5^2$ & $2\,1,2,2$ & $(2p_1)\,1,(2p_2),(2p_3)$ & $2p_2+2p_3-2p_1+1$ & $$ & $$   \\ \hline
$7_6^2$ & $.2$ & $.(2p_1)$ & $-2p_1+1$ & $$ & $$   \\ \hline
$8_4^2$ & $3\,2\,3$ & $(2p_1+1)\,(2p_2)\,(2p_3+1)$ & $-2p_1-2p_3-1$ & $$ & $$   \\ \hline
$8_5^2$ & $3\,1\,2\,2$ & $(2p_1+1)\,1\,(2p_2)\,(2p_3)$ & $-2p_1-1$ & $$ & $$   \\ \hline
$8_7^2$ & $2\,1\,2\,1\,2$ & $(2p_1)\,1\,(2p_2)\,1\,(2p_3)$ & $-2p_1-2p_3+1$ & $p_1+p_3$ & $p_1+p_3\ge p_2$   \\ \hline
$8_9^2$ & $2\,2,2,2$ & $(2p_1)\,(2p_2),(2p_3),(2p_4)$ & $-2p_3-2p_4+1$ & $$ & $$   \\ \hline
$8_{10}^2$ & $2\,1\,1,2,2$ & $(2p_1)\,1\,1,(2p_2),(2p_3)$ & $2p_3-1$ & $p_3$ & $p_1=p_2=1$   \\ \hline
$8_{11}^2$ & $3,2,2+$ & $(2p_1+1),(2p_2),(2p_3)+$ & $-2p_1+1$ & $$ & $$   \\ \hline
$8_{12}^2$ & $2\,1,2,2+$ & $(2p_1)\,1,(2p_2),(2p_3)+$ & $-2p_1+2p_2-1$ & $p_1-p_2+1$ & $p_2=1$   \\ \hline
\end{tabular}

\noindent \begin{tabular}{|c|c|c|c|c|c|} \hline
$K$ & $Con$ & $Fam$ & $\sigma $ & $u$ & $Cond$   \\ \hline
$8_{13}^2$ & $.2\,1$ & $.(2p_1)\,1$ & $-1$ & $$ & $$   \\ \hline
$8_{14}^2$ & $.2:2$ & $.(2p_1):(2p_2)$ & $1$ & $$ & $$   \\ \hline
$9_6^2$ & $3\,3\,1\,2$ & $(2p_1+1)\,(2p_2+1)\,1\,(2p_3)$ & $-2p_1+2p_3-1$ & $p_3-p_1$ & $p_3\ge p_1+p_2$   \\ \hline
$9_7^2$ & $3\,2\,1\,1\,2$ & $(2p_1+1)\,(2p_2)\,1\,1\,(2p_3)$ & $-2p_1+1$ & $p_1+1$ & $p_2=1,p_3\le p_1+2$    \\ \hline
$9_8^2$ & $3\,1\,3\,2$ & $(2p_1+1)\,1\,(2p_2+1)\,(2p_3)$ & $-2p_1+2p_3-1$ & $p_1-1$ & $p_1\ge p_2+p_3$   \\ \hline
$9_9^2$ & $3\,1\,1\,1\,3$ & $(2P-1+1)\,1\,1\,1\,(2p_2+1)$ & $-2p_1+1$ & $p_1$ & $p_1-p_2\ge 1$   \\ \hline
$9_{11}^2$ & $2\,2\,2\,1\,2$ & $(2p_1)\,(2p_2)\,(2p_3)\,1\,(2p_4)$ & $-2p_1-2p_3+1$ & $p_1+p_3$ & $p_3\ge p_4$   \\ \hline
$9_{12}^2$ & $2\,2\,1\,1\,1\,2$ & $(2p_1)\,(2p_2)\,1\,1\,1\,(2p_3)$ & $-2p_1+2p_3-1$ & $p_3-p_1$ & $p_3>p_1+p_2$   \\
$$ & $$ & $$ & $$ & $p_1-p_3+1$ & $p_2=1, p_1>p_3$   \\ \hline
$9_{15}^2$ & $3\,2,2,2$ & $(2p_1+1)\,(2p_2),(2p_3),(2p_4)$ & $-2p_2-2p_4+1$ & $$ & $$   \\ \hline
$9_{16}^2$ & $3\,1\,1,2,2$ & $(2p_1+1)\,1\,1,(2p_2),(2p_3)$ & $2p_2+2p_3-3$ & $$ & $$   \\ \hline
$9_{17}$ & $2\,3,2,2$ & $(2p_1)\,(2p_2)\,1,(2p_3),(2p_4)$ & $-2p_2+2p_3+2p_4-1$ & $$ & $$   \\ \hline
$9_{18}^2$ & $2\,2\,1,2,2$ & $(2p_1)\,(2p_2)\,1,(2p_3),(2p_4)$ & $-2p_1-2p_4+1$ & $p_1+p_4$ & $p_3=1$   \\ \hline
$9_{21}^2$ & $3\,1,3,2$ & $(2p_1+1)\,1,(2p_2+1),(2p_3)$ & $-2p_1+2p_2+2p_3-1$ & $p_2+p_3-p_1$ & $p_3>p_1$   \\ \hline
$9_{22}^2$ & $3\,1,2\,1,2$ & $(2p_1+1)\,1,(2p_2)\,1,(2p_3)$ & $-2p_1-2p_2-1$ & $$ & $$   \\ \hline
$9_{23}^2$ & $3,3,2\,1$ & $(2p_1+1),(2p_2+1)\,(2p_3)\,1$ & $-2p_1-2p_2+2p_3-1$ & $p_1+p_2-p_3+1$ & $p_1\ge p_3$ or $p_2\ge p_3$  \\ \hline
$9_{24}^2$ & $2\,1,2\,1,2\,1$ & $(2p_1)\,1,(2p_2)\,1,(2p_3)\,1$ & $-2p_1-2p_2-2p_3-1$ & $p_1+p_2+p_3$ & $$   \\ \hline
$9_{25}^2$ & $2\,2,2,2+$ & $(2p_1)\,(2p_2),(2p_3),(2p_4)+$ & $-2p_2+1$ & $$ & $$   \\ \hline
$9_{26}^2$ & $2\,1\,1,2,2+$ & $(2p_1)\,1\,1,(2p_2),(2p_3)+$ & $-1$ & $$ & $$   \\ \hline
$9_{27}^2$ & $3,2,2+2$ & $(2p_1+1),(2p_2),(2p_3)+(2p_4)$ & $-2p_1-2p_3+1$ & $$ & $$   \\ \hline
$9_{28}^2$ & $2\,1,2,2,2+2$ & $(2p_1)\,1,(2p_2),(2p_3)+(2p_4)$ & $-2p_1+2p_2+2p_3-1$ & $$ & $$   \\ \hline
$9_{29}^2$ & $(3,2)\,(2,2)$ & $(2p_1+1),(2p_2))\,((2p_3),(2p_4))$ & $2p_1+2p_2-2p_3-2p_4-1$ & $$ & $$   \\ \hline
$9_{30}^2$ & $(2\,1,2)\,(2,2)$ & $((2p_1)\,1,(2p_2))\,((2p_3),(2p_4))$ & $-2p_1-2p_3-2p_4+1$ & $p_1+p_3+p_4$ & $$   \\ \hline
$9_{32}^2$ & $.3\,1$ & $.(2p_1+1)\,1$ & $-1$ & $$ & $$   \\ \hline
$9_{33}^2$ & $.2\,2$ & $.(2p_1)\,(2p_2)$ & $-2p_1+1$ & $$ & $$   \\ \hline
$9_{34}^2$ & $.3.2$ & $.(2p_1+1).(2p_2)$ & $-2p_1+2p_2-1$ & $p_1$ & $p_1>1,p_2=1$   \\ \hline
$9_{35}^2$ & $.3.2\,0$ & $.(2p_1+1).(2p_2)\,0$ & $-2p_1-2p_2+1$ & $p_1+p_2$ & $$   \\ \hline
$9_{36}^2$ & $.3:2$ & $.(2p_1+1):(2p_2)$ & $-2p_1+2p_2-1$ & $p_1-p_2+1$ & $p_1>1,p_2=1$   \\ \hline
$9_{37}^2$ & $.3:2\,0$ & $.(2p_1+1):(2p_2)\,0$ & $-2p_1+2p_2-1$ & $$ & $$   \\ \hline
$9_{38}^2$ & $.2\,1:2\,0$ & $.(2p_1)\,1:(2p_2)\,0$ & $-2p_1+3$ & $$ & $$   \\ \hline
$9_{39}^2$ & $.2.2.2\,0$ & $.(2p_1).(2p_2).(2p_3)\,0$ & $-2p_1+2p_2+2p_3-1$ & $p_2+p_3-p_1$ & $p_3>p_1$   \\ \hline
$9_{40}^2$ & $2:2:2$ & $(2p_1):(2p_2):(2p_3)$ & $-2p_1-2p_2-2p_3+3$ & $$ & $$   \\ \hline
$9_{41}^2$ & $2:2\,0:2\,0$ & $(2p_1):(2p_2):(2p_3)\,0$ & $-2p_1+2p_2+2p_3-1$ & $p_2+p_3-p_1$ & $p_2>p_1$ or $p_3>p_1$   \\ \hline
$9_{42}^2$ & $8^*2$ & $8^*(2p_1)$ & $-2p_1+1$ & $p_1$ & $p_1\ge 2$   \\ \hline
$6_1^3$ & $2,2,2$ & $(2p_1),(2p_2),(2p_3)$ & $-2p_1+2p_2-2p_3$ & $$ & $$   \\ \hline
$7_1^3$ & $2,2,2+$ & $(2p_1),(2p_2),(2p_3)+$ & $-2p_1+2$ & $$ & $$   \\ \hline
$8_3^2$ & $3\,1,2,2$ & $(2p_1+1)\,1,(2p_2),(2p_3)$ & $-2p_1+2p_2+2p_3-2$ & $$ & $$   \\ \hline
$8_3^3$ & $2,2,2+2$ & $(2p_1),(2p_2),(2p_3)+(2p_4)$ & $-2p_1-2p_3+2$ & $$ & $$   \\ \hline
$8_4^3$ & $(2,2)\,(2,2)$ & $((2p_1),(2p_2))\,((2p_3),(2p_4))$ & $-2p_1-2p_2+2p_3+2p_4$ & $$ & $$   \\ \hline
$8_5^3$ & $.3$ & $.(2p_1+1)$ & $-2p_1$ & $$ & $$   \\ \hline
$8_6^3$ & $.2:2\,0$ & $.(2p_1):(2p_2)\,0$ & $-2p_1+2p_2$ & $$ & $$   \\ \hline

\end{tabular}

\end{landscape}

\noindent \begin{tabular}{|c|c|c|c|c|c|} \hline
$9_1^3$ & $2\,1\,2,2,2$ & $(2p_1)\,1\,(2p_2),(2p_3),(2p_4)$ & $-2p_1-2p_3-2p_4+2$ & $$ & $$   \\ \hline
$9_2^3$ & $2\,1\,1\,1,2,2$ & $(2p_1)\,1\,1\,1,(2p_2),(2p_3)$ & $-2p_1+2p_3$ & $$ & $$   \\ \hline
$9_3^3$ & $3,2,2,2$ & $(2p_1+1),(2p_2),(2p_3),(2p_4)$ & $-2p_1-2p_4+2$ & $$ & $$   \\ \hline
$9_4^3$ & $2\,1,2,2,2$ & $(2p_1)\,1,(2p_2),(2p_3),(2p_4)$ & $-2p_1+2p_2+2p_4-2$ & $$ & $$   \\ \hline
$9_6^3$ & $3\,1,2,2$ & $(2p_1+1)\,1,(2p_2),(2p_3)+$ & $-2p_1+2p_2-2$ & $$ & $$   \\ \hline
$9_7^3$ & $2,2,2+3$ & $(2p_1),(2p_2),(2p_3)+(2p_4+1)$ & $-2p_1+2$ & $$ & $$   \\ \hline
$9_8^3$ & $(2,2+)\,(2,2)$ & $((2p_1),(2p_2)+)\,((2p_3),(2p_4))$ & $-2p_1+2p_3+2p_4$ & $$ & $$   \\ \hline
$9_9^3$ & $(2,2)\,1\,(2,2)$ & $((2p_1),(2p_2)\,1\,((2p_3),(2p_4))$ & $-2p_1-2p_4+2$ & $$ & $$   \\ \hline
$9_{10}^3$ & $.2\,1\,1$ & $.(2p_1)\,1\,1$ & $-2p_1$ & $$ & $$   \\ \hline
$9_{11}^3$ & $.2\,1:2$ & $(-2p_1)\,1:2$ & $-2p_1+2$ & $$ & $$   \\ \hline
$9_{12}^3$ & $.(2,2)$ & $.((2p_1),(2p_2))$ & $-2p_1+2$ & $$ & $$   \\ \hline
$8_1^4$ & $2,2,2,2$ & $(2p_1),(2p_2),(2p_3),(2p_4)$ & $-2p_1-2p_4+3$ & $$ & $$   \\ \hline

\end{tabular}

\bigskip

\normalsize

\section{Ascending numbers of alternating knot families}

Our next goal was to compute ascending numbers of alternating knots belonging to some families with known unknotting numbers and to find their based oriented diagrams showing the ascending number. These results for the families beginning with knots with $n\le 8$ crossings are described in Theorems 5.1-5.4, for the families beginning with knots with $n=9$ crossings in Theorems 5.5-5.9, and for the families beginning with knots with $n=10$ crossings in Theorems 5.10-5.15.

For all these families, except for the first and the last family, ascending numbers are computable only from non-minimal diagrams. Hence, in Figs. 29-38 every family is represented by its minimal diagram (a), and non-minimal based oriented diagram (b) giving the corresponding ascending number.

\begin{theorem}
For knots $3_1=3$, $5_1=5$, $7_1=7$, $9_1=9$, $\ldots $ of the family $2p+1$ ($p\ge 1$), the minimal ascending number is $a(K)=u(K)=p$, and it is realized on the minimal diagrams (Fig. 29).
\end{theorem}

\begin{theorem}
For knots $7_3=4\,3$, $9_3=6\,3$, $\ldots $ of the family $(2p)\,3$ ($p\ge 2$) the minimal ascending number is $a(K)=u(K)=p$, and it is realized on the non-minimal diagrams of the form $((((1,(-1,(((1^{2p-2}),-1),-1))),1),-1),-1,-1)$ (Fig. 30b), where by $1^{2p-2}$ is denoted the sequence $1,\ldots ,1$ of the length $2p-2$.
\end{theorem}

\begin{theorem}
For knots $7_5=3\,2\,2$, $9_6=5\,2\,2$, $9_9=4\,2\,3$, $\ldots $ of the family $(2p+1)\,2\,(2q)$ ($p\ge 1$, $q\ge 1$), the minimal ascending number is $a(K)=u(K)=p+q$, and it is itself on the non-minimal diagrams of the form $(2p+1),-2\,1,(2q)$ (Fig. 31b).
\end{theorem}

\begin{theorem}
For knots $8_{15}=2\,1,2\,1,2$, $10_{49}=4\,1,2\,1,2$, $\ldots $ of the family $(2p)\,1,$ $(2q)\,1,2$ ($p\ge 1$, $q\ge 1$) minimal ascending number is $a(K)=u(K)=p+q$, and it is itself on the non-minimal diagrams of the form $(2p)\,1,(2q)\,1,-2,1$ (Fig. 32b).
\end{theorem}

For the knot family $2p+1$ ($p\ge 1$) the absolute value of the signature of a knot $2p+1$ ($p\ge 1$) is $2p$. For the knot family $(2p)\,3$ ($p\ge 2$) the absolute value of the signature is $2p$. For the family of knots $(2p+1)\,2\,(2q)$ ($p\ge 1$, $q\ge 1$) the absolute value of the signature is $2p+2q$. For the family of knots $(2p)\,1,(2q)\,1,2$ ($p\ge 1$, $q\ge 1$) the absolute value of the signature is $2p+2q$.

K.~Murasugi \cite{9} proved the lower bound for the unknotting number for knots $u(K)\geq  \frac{|\sigma _K|}{2}$. For the family of knots $2p+1$ ($p\ge 1$),
from Murasugi's Theorem 1.1, half of the signature is $p$, for the family of knots $(2p)\,3$ ($p\ge 2$) from Theorem 1.2 half of the signature is $p$, for the family of knots $(2p+1)\,2\,(2q)$ ($p\ge 1$, $q\ge 1$) from Theorem 1.3 it is $p+q$, and for the family of knots $(2p+1)\,2\,(2q)$ ($p\ge 1$, $q\ge 1$) from Theorem 1.4 half of the signature is $p+q$. The proof of the Theorems 1.1-1.4 follows from the fact that in each of the families the unknotting number is equal to a half of the absolute value of the signature, and it is realized on the minimal diagrams of the knots belonging to the families from Theorems 1.1-1.4, respectively. Hence, unknotting number for these families is equal to the minimal diagram unknotting number. For all these families, we effectively constructed corresponding diagrams giving minimal ascending number equal to the unknotting number, so it follows that for each of the knots from these families $a(K)=u(K)=\frac{|\sigma _K|}{2}$.

We provide similar theorems for certain families beginning with nine and ten crossing knots. Their proof is analogous to the proof of Theorems 5.1-5.4. For rational knot family $5\,(2p)$ ($p\ge 2$) from Theorem 5.5 the absolute value of the signature is $2p$; for rational knot family $(2p+1)\,4\,(2q)$ ($p\ge 1$, $q\ge 1$) from Theorem 5.6 the absolute value of the signature is $2p+2q$; for knot family $(2p+1)\,4\,(2q)$ ($p\ge 1$, $q\ge 1$) from Theorem 5.7 the absolute value of the signature is $2p+4$; for rational knot family $3\,2\,2\,(2p)$ ($p\ge 1$) from Theorem 5.8  the absolute value of the signature is $2p+2$; for rational knot family $(2p)\,2\,1\,2\,(2q)$ ($p\ge 1$, $q\ge 1$) from Theorem 5.9 the absolute value of the signature is $2p+2q$; for knot family $(2p)\,3,(2q)\,1,2$ ($p\ge 1$, $q\ge 1$) from Theorem 5.10  the absolute value of the signature is $2p+2q$; for knot family $(2p)\,3,(2q)\,1,2$ ($p\ge 1$, $q\ge 1$) from Theorem 5.11 the absolute value of the signature is $2p+4$. By similar arguments, applied to non-alternating knots we can conclude that for knot family $(2p)\,2\,1,(2q+1),-2$ ($p\ge 1$, $q\ge 1$) from Theorem 5.12 the absolute value of the signature is $2p+2q$, and for knot family $-3\,0:(2p)\,0:2\,0$ ($p\ge 1$) from Theorem 5.13 the absolute value of the signature is $2p+2$.

\begin{theorem}
For knots $9_4=5\,4$, $\ldots $ of the family $5\,(2p)$ ($p\ge 2$), the minimal ascending number is $a(K)=u(K)=p$, and it is realized on the non-minimal diagrams of the form $((1,(-1,(((((1,(1,(-1,-1))),1),-1),-1),-1))),1^{2p-2})$ (Fig. 33b), where by $1^{2p-2}$ is denoted the sequence $1,\ldots ,1$ of the length $2p-2$.
\end{theorem}

\begin{theorem}
For knots $9_7=3\,4\,2$, $\ldots $ of the family $(2p+1)\,4\,(2q)$ ($p\ge 1$, $q\ge 1$) the minimal ascending number is $a(K)=u(K)=p+q$, and it is realized on the non-minimal diagrams of the form $(((-1,(1,((((-1,((1^{2p}),1)),-1),1),1))),-1),$ $(-1)^{2q})$ (Fig. 34b), where by $1^{2p-2}$ and $(-1)^{2q-2}$ are denoted the sequence $1,\ldots ,1$ of the length $2p-2$, and $-1,\ldots ,-1$ of the length $2q-2$, respectively.
\end{theorem}

\begin{theorem}
For knots $9_{16}=3,3,2+$, $\ldots $ of the family $3,3,(2p)+$ ($p\ge 1$) the minimal ascending number is $a(K)=u(K)=p+2$, and it is realized on the non-minimal diagrams of the form $-(1,1)\,1\,1,-(1,1)\,1\,1,((-1)^{2p+1})\,1$ (Fig. 35b), where by $1^{2p+1}$ is denoted the sequence $1,\ldots ,1$ of the length $2p+1$.
\end{theorem}

\begin{theorem}
For knots $9_{18}=3\,2\,2\,2$, $\ldots $ of the family $3\,2\,2\,(2p)$ ($p\ge 1$), the minimal ascending number is $a(K)=u(K)=p+1$, and it is realized on non-minimal diagrams of the form $(((((1,(-1,(((-1,(1,1)),-1),-1))),1),-1),-1),(-1)^{2p})$ (Fig. 36b), where by $1^{2p}$ is denoted the sequence $-1,\ldots ,-1$ of the length $2p$.
\end{theorem}

\begin{theorem}
For knots $9_{23}=2\,2\,1\,2\,2$, $\ldots $ of the family $(2p)\,2\,1\,2\,(2q)$ ($p\ge 1$, $q\ge 1$), the minimal ascending number is $a(K)=u(K)=p+q$, and it is realized on the non-minimal diagrams of the form $((((((((1,(1,(1^{2p}))),-1),-1),1),1),-1),-1),$ $(-1)^{2q})$ (Fig. 37b),  where by $1^{2p}$ and $(-1)^{2q}$ are denoted the sequence $1,\ldots ,1$ of the length $2p$, and $-1,\ldots ,-1$ of the length $2q$, respectively.
\end{theorem}

\begin{theorem}
For knots $10_{50}=2\,3,2\,1,2$, $\ldots $ of the family $(2p)\,3,(2q)\,1,2$ ($p\ge 1$, $q\ge 1$), the minimal ascending number is $a(K)=u(K)=p+q$, and it is realized on the non-minimal diagrams of the form $(((-1,(-1,-(1^{2p}))),1),1),((((-(1^{2q}),1),-1),1),$ $1),-(1,1),1$ (Fig. 38b),  where by $1^{2p}$ and $1^{2q}$ are denoted the sequence $1,\ldots ,1$ of the length $2p$, and $1,\ldots ,1$ of the length $2q$, respectively.
\end{theorem}

\begin{theorem}
For knots $10_{80}=(3,2)\,(2\,1,2)$, $\ldots $ of the family $(3,2)\,((2p)\,1,2)$ ($p\ge 1$), the minimal ascending number is $a(K)=u(K)=p+2$, and it is realized on the non-minimal diagrams of the form $(((1,-(1,1)),1),-(1,1),1)$ $(((((-(1^{2p}),1),-1),1),1),$ $-(1,1),1)$ (Fig. 39b),  where by $1^{2p}$ is denoted the sequence $1,\ldots ,1$ of the length $2p$.
\end{theorem}

\begin{figure}[th]
\centerline{\psfig{file=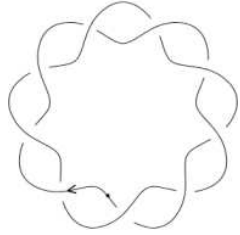,width=1.2in}}
\vspace*{8pt}
\caption{Family $2p+1$ ($p\ge 1$) and its corresponding minimal based oriented diagram giving the ascending number $p$.  \label{fig29}}
\end{figure}

\begin{figure}[th]
\centerline{\psfig{file=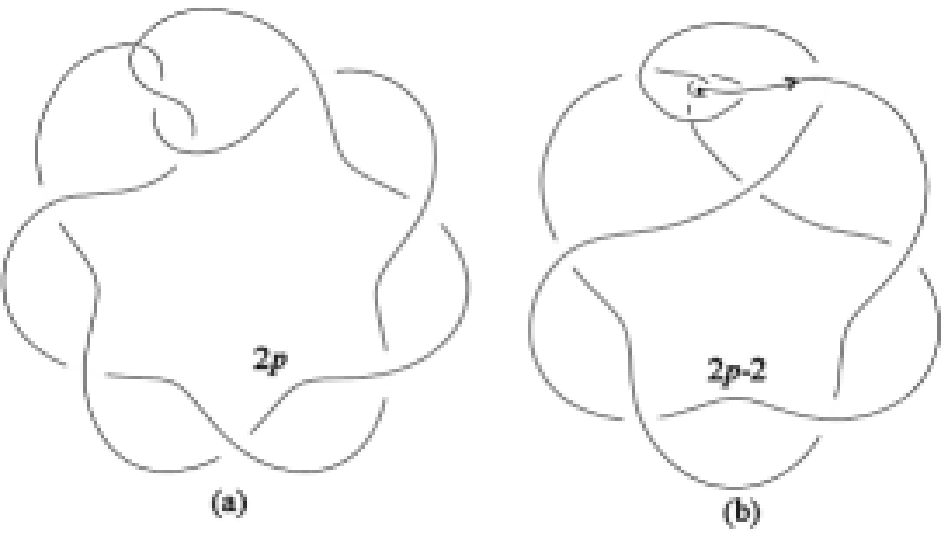,width=2.4in}}
\vspace*{8pt}
\caption{(a) Family $(2p)\,3$ ($p\ge 2$) with the ascending number $p$ (b). \label{fig30}}
\end{figure}

\begin{figure}[th]
\centerline{\psfig{file=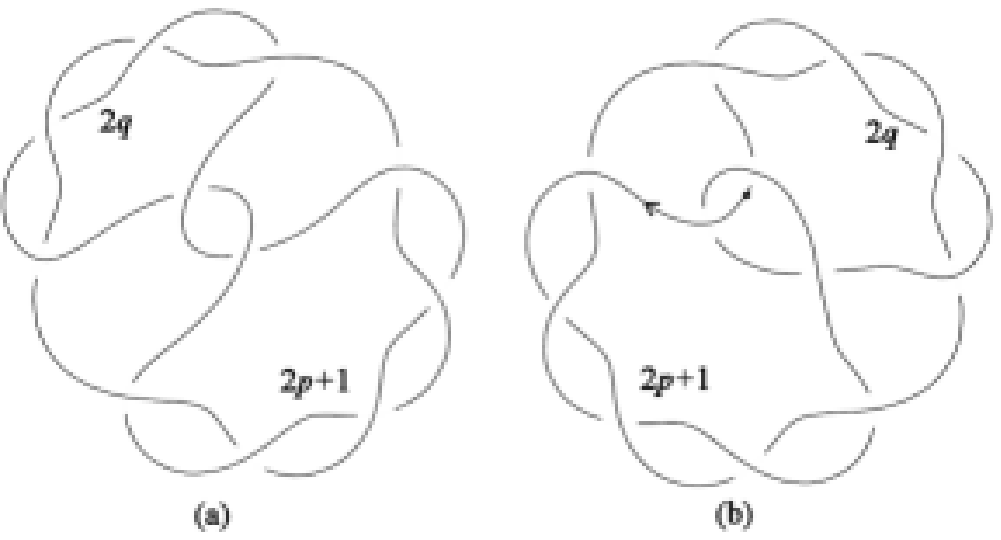,width=2.4in}}
\vspace*{8pt}
\caption{(a) Family $(2p+1)\,2\,(2q)$ ($p\ge 1$, $q\ge 1$) with the ascending number $p+q$ (b) . \label{fig31}}
\end{figure}

\begin{figure}[th]
\centerline{\psfig{file=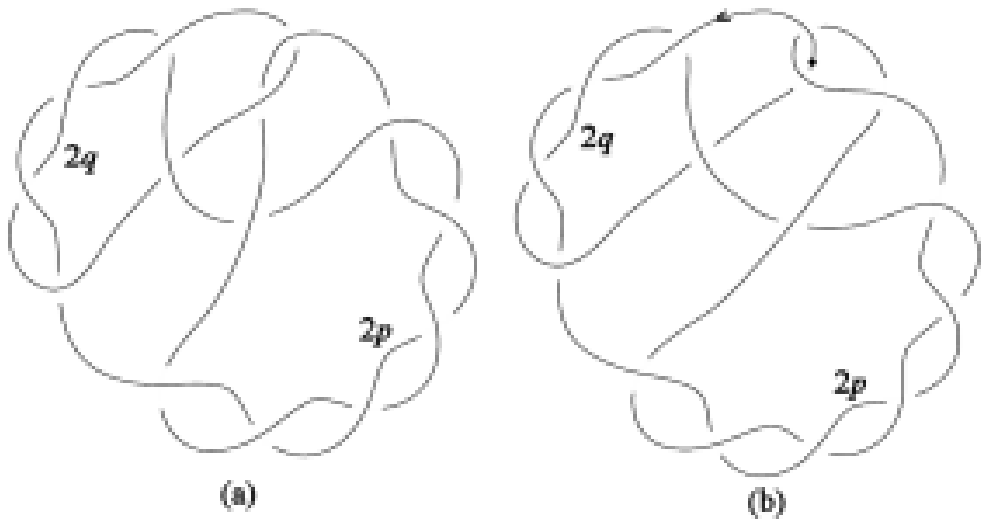,width=2.4in}}
\vspace*{8pt}
\caption{(a) Family $(2p)\,1,(2q)\,1,2$ ($p\ge 1$, $q\ge 1$)  with the ascending number $p+q$  (b) . \label{fig32}}
\end{figure}

\begin{figure}[th]
\centerline{\psfig{file=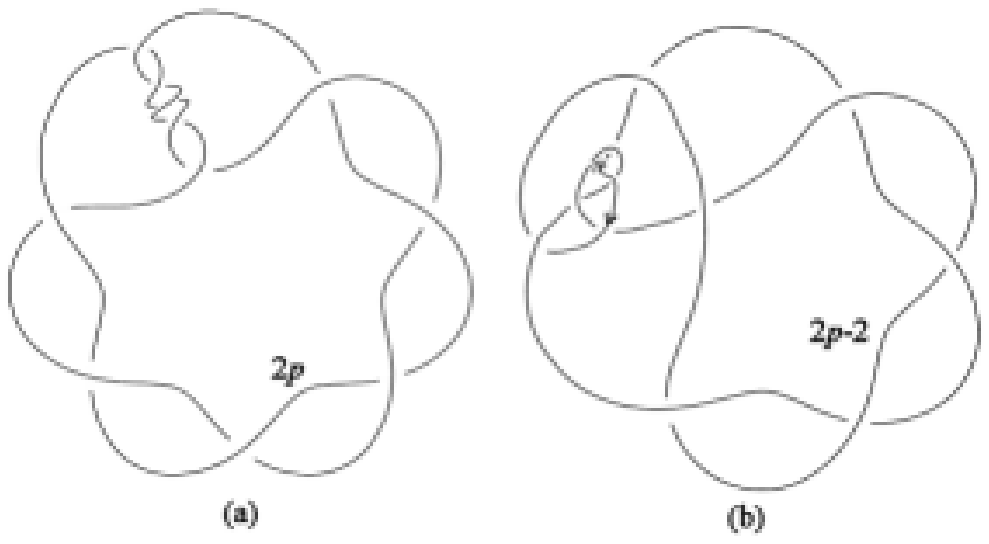,width=2.4in}}
\vspace*{8pt}
\caption{(a) Family $5\,(2p)$ ($p\ge 2$) with the ascending number $p$ (b). \label{fig33}}
\end{figure}

\begin{figure}[th]
\centerline{\psfig{file=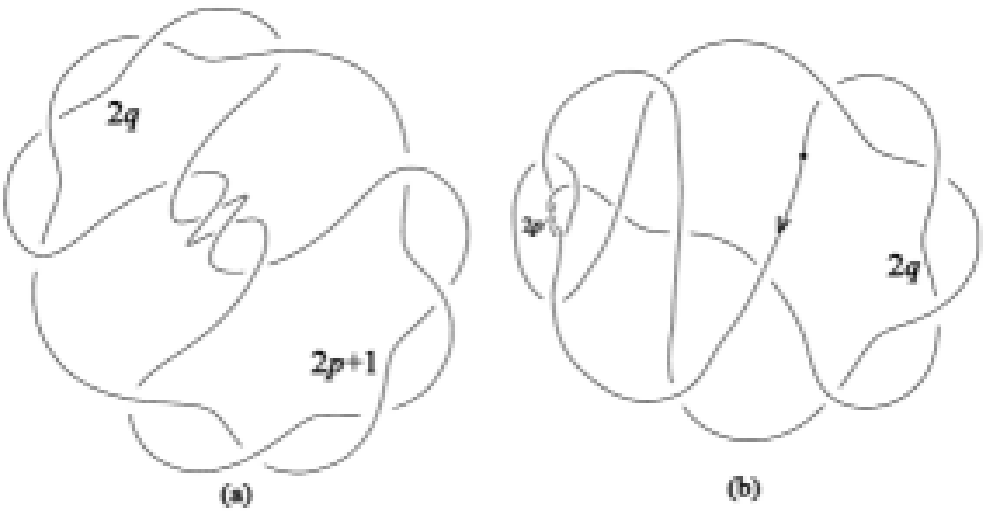,width=2.4in}}
\vspace*{8pt}
\caption{(a) Family $(2p+1)\,4\,(2q)$ ($p\ge 1$, $q\ge 1$)  with the ascending number $p+q$ (b). \label{fig34}}
\end{figure}

\begin{figure}[th]
\centerline{\psfig{file=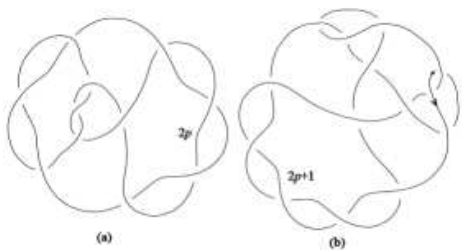,width=2.4in}}
\vspace*{8pt}
\caption{(a) Family $3,3,(2p)+$ ($p\ge 1$)  with the ascending number $p+2$ (b). \label{fig35}}
\end{figure}

\begin{figure}[th]
\centerline{\psfig{file=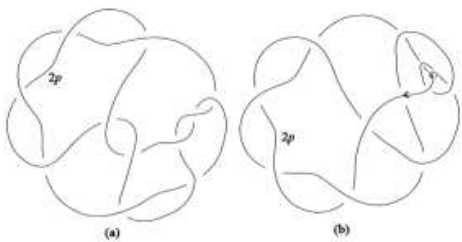,width=2.4in}}
\vspace*{8pt}
\caption{(a) Family  $3\,2\,2\,(2p)$ ($p\ge 1$) with the ascending number $p+1$ (b). \label{fig36}}
\end{figure}

\begin{figure}[th]
\centerline{\psfig{file=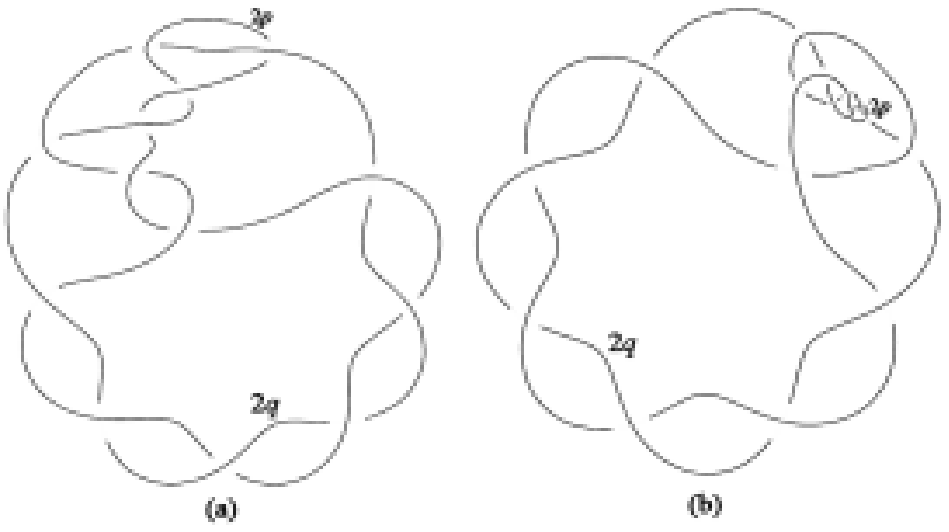,width=2.8in}}
\vspace*{8pt}
\caption{(a) Family $(2p)\,2\,1\,2\,(2q)$ ($p\ge 1$, $q\ge 1$) with the ascending number $p+q$ (b). \label{fig37}}
\end{figure}

\begin{figure}[th]
\centerline{\psfig{file=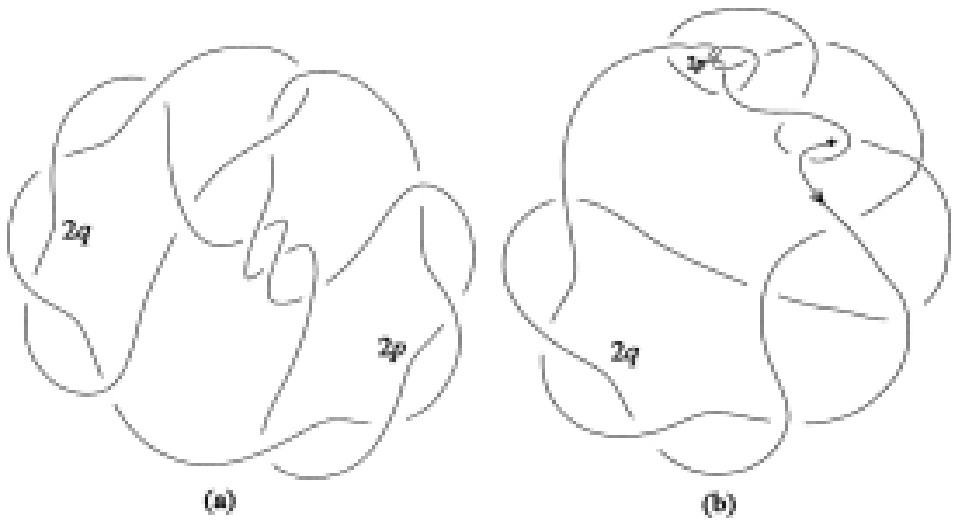,width=2.8in}}
\vspace*{8pt}
\caption{(a) Family $(2p)\,3,(2q)\,1,2$ ($p\ge 1$, $q\ge 1$) with the ascending number $p+q$ (b). \label{fig38}}
\end{figure}

\begin{figure}[th]
\centerline{\psfig{file=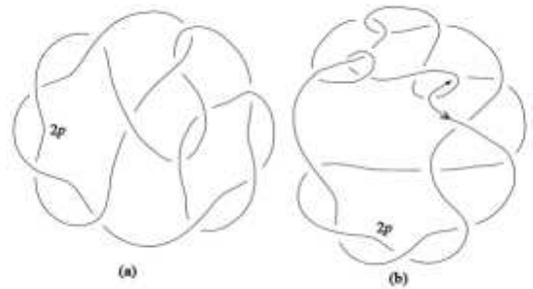,width=2.8in}}
\vspace*{8pt}
\caption{(a) Family  $(3,2)\,((2p)\,1,2)$ ($p\ge 1$)  with the ascending number $p+2$ (b). \label{fig39}}
\end{figure}

\end{document}